\documentclass[10pt]{article}
\usepackage{amsmath}
\usepackage{amssymb}
\def\D{\displaystyle}

\begin{document}

\def\D{\displaystyle}
\newtheorem{theorem}{Theorem}[section]
\newtheorem{definition}[theorem]{Definition}
\newtheorem{lemma}[theorem]{Lemma}
\newtheorem{proposition}[theorem]{Proposition}
\newtheorem{remark}[theorem]{Remark}
\renewcommand{\theequation}{\thesection. \arabic{equation}}
\newtheorem{example}[theorem]{Example}
\begin{center}
\large{\bf Multivariate Difference-Differential Dimension Polynomials}

\medskip

\normalsize Alexander Levin

Department of Mathematics, The Catholic University of America

Washington, D. C. 20064\\ E-mail: levin@cua.edu
\end{center}

\begin{abstract}
In this paper we generalize the Ritt-Kolchin method of characteristic sets and the
classical Gr\"obner basis technique to prove the existence and obtain methods of computation
of multivariate difference-differential dimension polynomials associated
with a finitely generated difference-differential field extension. We also give an
interpretation of such polynomials in the spirit of the A. Einstein's concept of strength
of a system of PDEs and determine their invariants, that is, characteristics of a
finitely generated difference-differential field extension carried by every
its dimension polynomial.
\end{abstract}

\section{Introduction}
The role of Hilbert polynomials in commutative and homological algebra as
well as in algebraic geometry and combinatorics is well known. A similar role
in differential algebra is played by differential dimension polynomials, which
describe in exact terms the freedom degree of a dynamic system
as well as the number of arbitrary constants in the general solution of a
system of partial algebraic differential equations.

The notion of differential dimension polynomial was introduced by E. Kolchin
in 1964 \cite{K1} who proved the following fundamental result.

\begin{theorem} Let $K$ be a differential field (\,$Char\,K = 0$), that is, a field considered together
with the action of a set $\Delta = \{\delta_{1},\dots, \delta_{m}\}$ of
mutually commuting derivations of $K$ into itself. Let $\Theta$ denote the free commutative
semigroup of all power products of the form $\theta = \delta_{1}^{k_{1}}\dots
\delta_{m}^{k_{m}}$ ($k_{i}\geq 0$), let $ord\,\theta = \sum_{i=1}^{m}k_{i}$, and for any
$r\geq 0$, let  $\Theta(r) = \{\theta\in \Theta\,|\,
ord\,\theta\leq r\}$. Furthermore, let $L = K\langle \eta_{1},\dots,\eta_{n}\rangle$ be a differential
field extension of $K$ generated by a finite set $\eta =
\{\eta_{1}, \dots , \eta_{n}\}$. \,(As a field, $L = K(\{\theta
\eta_{j} | \theta \in \Theta, 1\leq j\leq n\})$.\,)

Then there exists a polynomial
$\omega_{\eta|K}(t)\in {\bf Q}[t]$ such that

 {\em (i)}\, $\omega_{\eta|K}(r) = trdeg_{K}K(\{\theta \eta_{j} | \theta
\in \Theta(r), \,1\leq j\leq n\})$ for all sufficiently large $r\in {\bf
Z}$;

{\em (ii)}\, $\deg \omega_{\eta|K} \leq m$ and $\omega_{\eta|K}(t)$ can be
written as $\omega_{\eta|K}(t) = \D\sum_{i=0}^{m}a_{i}\D{t+i\choose i}$

\noindent where $a_{0},\dots, a_{m}\in {\bf Z}$;

{\em (iii)}\, $d =  \deg \omega_{\eta|K}$,\, $a_{m}$ and $a_{d}$ do not
depend on the choice of the system of $\Delta$-generators $\eta$
of the extension $L/K$ (clearly, $a_{d}\neq a_{m}$ iff $d < m$,
that is $a_{m}=0$).  Moreover, $a_{m}$ is equal to the
differential transcendence degree of $L$ over $K$, that is, to the
maximal number of elements $\xi_{1},\dots,\xi_{k}\in L$ such that
the set $\{\theta \xi_{i} | \theta \in \Theta, 1\leq i\leq k\}$ is
algebraically independent over $K$.
\end{theorem}

In 1980 A. Mikhalev and E. Pankrat'ev \cite{MP} showed that a system
of algebraic differential equations can be characterized by certain differential dimension polynomial, which
expresses the strength of the system in the sense of A. Einstein. The concept
of strength, which is an important characteristic of a system of PDEs governing
a physical field, was described by A. Einstein as follows
(see  \cite{E}): ''... the system of equations is to be chosen so that the
field quantities are determined as strongly as possible. In order to apply this
principle, we propose a method which gives a measure of strength of
an equation system. We expand the field variables, in the
neighborhood of a point $\mathcal{P}$, into a Taylor series (which
presupposes the analytic character of the field); the coefficients
of these series, which are the derivatives of the field variables at
$\mathcal{P}$, fall into sets according to the degree of
differentiation. In every such degree there appear, for the first
time, a set of coefficients which would be free for arbitrary choice
if it were not that the field must satisfy a system of differential
equations.  Through this system of differential equations (and its
derivatives with respect to the coordinates) the number of
coefficients is restricted, so that in each degree a smaller number
of coefficients is left free for arbitrary choice. The set of
numbers of ''free'' coefficients for all degrees of differentiation
is then a measure of the ''weakness'' of the system of equations,
and through this, also of its ''strength''.\,''

In this paper we generalize the Ritt-Kolchin method of characteristic sets to
the case of difference-differential polynomials and apply this method to prove the existence and
find invariants of multivariate dimension polynomials associated with a fixed partition of the
basic sets of derivations and automorphisms. We show that one can assign such a multivariate
polynomial to a system of partial algebraic difference-differential equations, and this polynomial
expresses the strength of the system in the sense of A. Einstein. We also find new invariants of
a finitely generated difference-differential field extension carried by multivariate dimension polynomials,
that is, characteristics of the extension, which do not depend on the choice of the system of its
generators.

\section{Preliminaries}

In this section we present some basic concepts and results used in the rest of the paper.

\smallskip

Throughout the paper, ${\bf N}, {\bf Z}$, ${\bf Q}$, and ${\bf R}$ denote the sets
of all non-negative integers, integers, rational numbers, and real numbers, respectively.
As usual, ${\bf Q}[t]$ denotes the ring of polynomials in one variable $t$
with rational coefficients. By a ring we always mean an associative ring with a
unity. Every ring homomorphism is unitary (maps unity onto unity), every subring
of a ring contains the unity of the ring. Unless otherwise indicated, by a module
over a ring $R$ we always mean a unitary left $R$-module.

\indent By a {\em difference-differential ring} we mean a commutative
ring $R$ together with finite sets $\Delta = \{\delta_{1},\dots,
\delta_{m}\}$ and $\sigma = \{\alpha_{1},\dots, \alpha_{n}\}$ of derivations
and automorphisms of $R$, respectively, such that any two mappings of the set
$\Delta\bigcup\sigma$ commute.  The set $\Delta\bigcup\sigma$ is called the
{\em basic set\/} of the difference-differential ring $R$, which is also called
a $\Delta\bigcup\sigma$-ring. If $R$ is a field, it is called a {\em difference-differential
field} or a $\Delta\bigcup\sigma$-field. Furthermore, in what follows, we denote the set
$\{\alpha_{1},\dots, \alpha_{n}, \alpha^{-1}_{1},\dots, \alpha^{-1}_{n}\}$ by $\sigma^{\ast}$.

\medskip

If $R$ is a difference-differential ring with a basic set $\Delta\bigcup\sigma$ described above, then
$\Lambda$ will denote the free commutative semigroup of all power products of the form
$\lambda = \delta_{1}^{k_{1}}\dots \delta_{m}^{k_{m}}\alpha_{1}^{l_{1}}\dots \alpha_{n}^{l_{n}}$ where
$k_{i}\in {\bf N},\, l_{j}\in {\bf Z}$ ($1\leq i\leq m,\, 1\leq j\leq n$). For any such an element $\lambda$,
we set $\lambda_{\Delta} = \delta_{1}^{k_{1}}\dots \delta_{m}^{k_{m}}$,
$\lambda_{\sigma} = \alpha_{1}^{l_{1}}\dots \alpha_{n}^{l_{n}}$, and denote by $\Lambda_{\Delta}$
and $\Lambda_{\sigma}$ the commutative semigroup of power products $\delta_{1}^{k_{1}}\dots \delta_{m}^{k_{m}}$
and the commutative group of elements of the form $\alpha_{1}^{l_{1}}\dots \alpha_{n}^{l_{n}}$, respectively.
The {\em order} of $\lambda$ is defined as $ord\,\lambda = \sum_{i=1}^{m}k_{i} + \sum_{j=1}^{n}|l_{j}|$, and for
every $r\in {\bf N}$, we set $\Lambda(r) = \{\lambda\in \Lambda\,|\, ord\,\lambda\leq r\}$ ($r\in {\bf N}$).

\medskip

A subring (ideal) $R_{0}$ of a $\Delta$-$\sigma$-ring $R$ is said to be a
difference-differential  (or $\Delta$-$\sigma$-) subring
of $R$ (respectively, difference-differential (or $\Delta$-$\sigma$-) ideal of $R$)
if $R_{0}$ is closed with respect to the action of any operator of $\Delta\bigcup\sigma^{\ast}$.
In this case the restriction of a mapping from $\Delta\bigcup\sigma^{\ast}$ on $R_{0}$ is denoted by same symbol.
If a prime (maximal) ideal $P$ of $R$ is closed with respect to the action of $\Delta\bigcup\sigma^{\ast}$,
it is called a {\em prime} (respectively, {\em maximal}) difference-differential (or $\Delta$-$\sigma$-) {\em ideal \/} of $R$.

If $R$ is a $\Delta$-$\sigma$-field and $R_{0}$ a subfield of $R$
which is also a $\Delta$-$\sigma$-subring of $R$, then  $R_{0}$ is said to be a
$\Delta$-$\sigma$-subfield of $R$; $R$, in turn, is called a
difference-differential (or $\Delta$-$\sigma$-) field extension or a $\Delta$-$\sigma$-overfield of $R_{0}$. In this case we also say
that we have a $\Delta$-$\sigma$-field extension $R/R_{0}$.

\medskip

If $R$ is a $\Delta$-$\sigma$-ring and $\Sigma \subseteq R$, then the intersection
of all $\Delta$-$\sigma$-ideals of $R$ containing the set $\Sigma$ is, obviously, the
smallest $\Delta$-$\sigma$-ideal of $R$ containing $\Sigma$. This ideal is denoted by $[\Sigma]$.
(It is clear that $[\Sigma]$ is generated, as an ideal, by the set
$\{\lambda \xi | \xi \in \Sigma,\, \lambda\in \Lambda\}$). If the set $\Sigma$ is finite,
$\Sigma = \{\xi_{1},\dots, \xi_{q}\}$, we say that the $\Delta$-$\sigma$-ideal
$I = [\Sigma]$ is finitely generated (we write this as
$I = [\xi_{1},\dots, \xi_{q}]$) and call $\xi_{1},\dots, \xi_{q}$
differential or $\Delta$-$\sigma$-generators of $I$.

If $K_0$ is a $\Delta$-$\sigma$-subfield of a $\Delta$-$\sigma$-field
$K$ and $\Sigma \subseteq K$, then the intersection of all $\Delta$-$\sigma$-subfields
of $K$ containing $K_0$ and $\Sigma$ is the unique $\Delta$-$\sigma$-subfield of
$K$ containing $K_0$ and $\Sigma$ and contained in every $\Delta$-$\sigma$-subfield
of $K$ containing $K_0$ and $\Sigma$. It is denoted by
$K_{0}\langle \Sigma \rangle$.
If $K = K_{0}\langle \Sigma \rangle$ and the set $\Sigma$ is finite,
$\Sigma = \{\eta_{1},\dots,\eta_{s}\}$, then $K$ is said to be a finitely
generated $\Delta$-$\sigma$-extension of $K_{0}$ with the set of $\Delta$-$\sigma$-generators
$\{\eta_{1},\dots,\eta_{s}\}$. In this case we write
$K = K_{0}\langle \eta_{1},\dots,\eta_{s}\rangle$.  It is easy to see that
the field $K_{0}\langle \eta_{1},\dots,\eta_{s}\rangle$ coincides with the
field $K_0(\{\lambda \eta_{i} | \lambda \in \Lambda, 1\leq i\leq s\}$).

\medskip

Let $R$ and $S$ be two difference-differential rings with the same basic set
$\Delta\bigcup\sigma$, so that elements of the sets $\Delta$ and $\sigma$ act on each
of the rings as mutually commuting derivations and automorphisms, respectively,
and every two mapping of the set $\Delta\bigcup\sigma$ commute.
(More rigorously, we assume that there exist injective mappings of the sets
$\Delta$ and $\sigma$ into the sets of derivations and automorphisms of the rings
$R$ and $S$, respectively, such that the images of any two elements of
$\Delta\bigcup\sigma$ commute. For convenience we will denote the images of elements of
$\Delta\bigcup\sigma$ under these mappings by the same symbols $\delta_{1},\dots, \delta_{m},
\alpha_{1},\dots, \alpha_{n}$). A ring homomorphism
$\phi: R \longrightarrow S$ is called a {\em difference-differential\/} or
$\Delta$-$\sigma$-{\em homomorphism\/} if $\phi(\tau a) = \tau \phi(a)$ for any
$\tau \in \Delta\bigcup\sigma$, $a\in R$. The notions of {\em $\Delta$-$\sigma$-epimorphism,
$\Delta$-$\sigma$-monomorphism, $\Delta$-$\sigma$-automorphism\/}, etc. are defined naturally
(as the corresponding ring homomorphisms that are $\Delta$-$\sigma$-homomorphisms).

\medskip

If $K$ is a difference-differential ($\Delta$-$\sigma$-) field and
$Y =\{y_{1},\dots, y_{s}\}$ is a finite set of symbols, then one can
consider the countable set of symbols $\Lambda Y = \{\lambda y_{j}|\lambda \in \Lambda, 1\leq j\leq s\}$ and the
polynomial ring $R = K[ \{\lambda y_{j}|\lambda \in \Lambda, 1\leq j\leq s\}]$
in the set of indeterminates $\Lambda Y$ over the field $K$.
This polynomial ring is naturally viewed as a $\Delta$-$\sigma$-ring where
$\tau(\lambda y_{j}) = (\tau\lambda)y_{j}$ for any $\tau \in \Delta\bigcup \sigma$,
$\lambda \in \Lambda$, $1\leq j\leq s$, and the elements of $\Delta$ act on
the coefficients of the polynomials of $R$ as they act in the field $K$.
The ring $R$ is called a {\em ring of difference-differential\/} (or $\Delta$-$\sigma$-)
{\em polynomials\/} in the set of differential ($\Delta$-$\sigma$-) indeterminates
$y_{1},\dots, y_{s}$ over $K$. This ring is denoted by
$K\{y_{1},\dots, y_{s}\}$ and its elements are called difference-differential
(or $\Delta$-$\sigma$-) polynomials.

\bigskip

Let $L = K\langle \eta_{1},\dots,\eta_{s}\rangle$ be a
difference-differential field extension of $K$ generated by a finite
set $\eta =\{\eta_{1},\dots,\eta_{s}\}$.

\medskip

The following is a unified version of E. Kolchin's theorem on
differential dimension polynomial and the author's theorem on the
dimension polynomial of a difference field extension
(see \cite{Levin1} or \cite[Theorem 4.2.5]{Levin4}\,).

\begin{theorem}  With the above notation, there exists a
polynomial $\phi_{\eta|K}(t)\in {\bf Q}[t]$ such that

\medskip

{\em (i)}\, $\phi_{\eta|K}(r) = trdeg_{K}K(\{\lambda\eta_{j} | \lambda\in
\Lambda(r), 1\leq j\leq s\})$ for all sufficiently large $r\in {\bf
Z}$;

\medskip

{\em (ii)}\, $\deg \phi_{\eta|K} \leq m+n$ and $\phi_{\eta|K}(t)$ can be
written as \, $\phi_{\eta|K}(t) =
\D\sum_{i=0}^{m+n}a_{i}{t+i\choose i}$ where $a_{0},\dots, a_{m+n}\in
{\bf Z}$ and $2^{n}|a_{m+n}$\,\,.

\medskip

{\em (iii)}\, $d =  \deg \phi_{\eta|K}$,\, $a_{m+n}$ and $a_{d}$ do not
depend on the set of difference-differential generators $\eta$ of
$L/K$ ($a_{d}\neq a_{m+n}$ if and only if $d < m+n$). Moreover,
$\D\frac{a_{m+n}}{2^{n}}$ is equal to the difference-differential
transcendence degree of $L$ over $K$ (denoted by
$\Delta$-$\sigma$-$trdeg_{K}L$), that is, to the maximal number of
elements $\xi_{1},\dots,\xi_{k}\in L$ such that the family
$\{\lambda \xi_{i} | \lambda \in \Lambda, 1\leq i\leq k\}$ is
algebraically independent over $K$.

\end{theorem}

The polynomial whose existence is established by this theorem is called a
{\em difference-differential} (or $\Delta$-$\sigma$-) {\em dimension polynomial} of
the extension $L/K$ associated with the system of difference-differential generators $\eta$.

\section{Partitions of the basic set of derivations and the formulation of the main theorem}
\setcounter{equation}{0}

Let $K$ be a difference-differential field of zero characteristic with basic
sets $\Delta = \{\delta_{1},\dots, \delta_{m}\}$ and $\sigma = \{\alpha_{1},\dots, \alpha_{n}\}$
of derivations and automorphisms, respectively. Suppose that the set of derivations is
represented as the union of $p$ disjoint subsets ($p\geq 1$):
\begin{equation}\Delta = \Delta_{1}\bigcup \dots \bigcup \Delta_{p}\end{equation}

where
$$\Delta_{1} = \{\delta_{1},\dots, \delta_{m_{1}}\},\, \Delta_{2} = \{\delta_{m_{1}+1},\dots, \delta_{m_{1}+m_{2}}\}, \,\dots,$$
$$\Delta_{p} = \{\delta_{m_{1}+\dots + m_{p-1}+1},\dots, \delta_{m}\}.\,  \, (m_{1}+\dots + m_{p} = m).$$

In other words, we fix a partition of the set $\Delta$ into $p$ sets of derivations.

\medskip

If \,\,$\lambda  = \delta_{1}^{k_{1}}\dots \delta_{m}^{k_{m}}\alpha_{1}^{l_{1}}\dots \alpha_{n}^{l_{n}}\in\Lambda$
($k_{i}\in {\bf N}, \,\, l_{j}\in {\bf Z}$), then the order of $\lambda$ with respect to a set $\Delta_{i}$  ($1\leq i\leq p$)
is defined as $\D\sum_{\nu = m_{1}+\dots + m_{i-1}+1}^{m_{1}+\dots + m_{i}}k_{\nu}$; it is denoted by $ord_{i}\lambda$.
(If $i=1$, the last sum is replaced by $k_{1}+\dots + k_{m_{1}}$.)
The number $ord_{\sigma}\lambda = \D\sum_{j = 1}^{n}|l_{j}|$ is called the order of $\lambda$ with respect to $\sigma$.

\medskip

If $r_{1},\dots, r_{p+1}\in {\bf N}$, we set

\medskip

$\Lambda(r_{1},\dots, r_{p+1}) = \{\lambda\in\Lambda\,|\,ord_{i}\lambda \leq r_{i}$
for $i=1,\dots, p$ and $ord_{\sigma}\lambda \leq r_{p+1}\}.$

\medskip

In what follows, for any permutation $(j_{1},\dots, j_{p+1})$ of the set
$\{1,\dots, p+1\}$, $<_{j_{1},\dots, j_{p+1}}$ will denote the lexicographic order
on ${\bf N}^{p+1}$ such that

\noindent$(r_{1},\dots, r_{p+1})<_{j_{1},\dots, j_{p+1}}
(s_{1},\dots, s_{p+1})$ if and only if either $r_{j_{1}} <
s_{j_{1}}$ or there exists $k\in {\bf N}$, $1\leq k\leq p$,
such that $r_{j_{\nu}} = s_{j_{\nu}}$ for $\nu = 1,\dots, k$ and
$r_{j_{k+1}} < s_{j_{k+1}}$. Furthermore, if $\Sigma \subseteq {\bf N}^{p+1}$,
then $\Sigma'$ denotes the set

\medskip

\noindent$\{e\in \Sigma | e$ is a maximal element of $\Sigma$ with
respect to one of the $(p+1)!$ lexicographic orders
$<_{j_{1},\dots, j_{p+1}}\}$. For example, if $\Sigma = \{(3, 0, 2), (2, 1, 1), (0, 1, 4), (1, 0, 3),$

\smallskip

\noindent$(1, 1, 6), (3, 1, 0), (1, 2, 0)\} \subseteq {\bf N}^{3}$, then
$\Sigma' = \{(3, 0, 2), (3, 1, 0), (1, 1, 6), (1, 2, 0)\}$.

\begin{theorem} Let $L = K\langle \eta_{1},\dots,\eta_{s}\rangle$ be a
$\Delta$-$\sigma$-field extension generated by a set $\eta =
\{\eta_{1}, \dots , \eta_{s}\}$. Then there exists a polynomial
$\Phi_{\eta}(t_{1},\dots, t_{p+1})$ in $(p+1)$ variables $t_{1},\dots, t_{p+1}$
with rational coefficients such that

\bigskip

{\em (i)} \,$\Phi_{\eta}(r_{1},\dots, r_{p+q}) = trdeg_{K}K(\D\bigcup_{j=1}^{s} \Lambda(r_{1},\dots,
r_{p+1})\eta_{j})$

\smallskip

\noindent for all sufficiently large $(r_{1},\dots,
r_{p+1})\in {\bf N}^{p+1}$ (i. e., there exist  $s_{1},\dots, s_{p+1}\in {\bf N}$ such
that the last equality holds for all $(r_{1},\dots,
r_{p+1})\in {\bf N}^{p+1}$ with $r_{1}\geq s_{1}, \dots, r_{p+1}\geq
s_{p+1}$);

\bigskip

{\em (ii)} \, $deg_{t_{i}}\Phi_{\eta} \leq m_{i}$ ($1\leq i\leq p$)\,\,\, and
$deg_{t_{p+1}}\Phi_{\eta} \leq n$, so that

\medskip

\noindent $deg\,\Phi_{\eta}\leq m+n$ and
$\Phi_{\eta}(t_{1},\dots, t_{p+1})$ can be represented as

\medskip

$\Phi_{\eta}(t_{1},\dots, t_{p+1}) = \D\sum_{i_{1}=0}^{m_{1}}\dots
\D\sum_{i_{p}=0}^{m_{p}}\D\sum_{i_{p+1}=0}^{n}a_{i_{1}\dots i_{p+1}}
{t_{1}+i_{1}\choose i_{1}}\dots {t_{p+1}+i_{p+1}\choose i_{p+1}}$

\bigskip

\noindent where $a_{i_{1}\dots i_{p+1}}\in {\bf Z}$ and $2^{n}\,|\,a_{m_{1}\dots m_{p}n}$.

\bigskip

{\em (iii)} \,Let $E_{\eta} = \{(i_{1},\dots, i_{p+1})\in {\bf N}^{p+1}\,|\,
0\leq i_{k}\leq m_{k}$ for $k=1,\dots, p$, $0\leq i_{p+1}\leq n$,
and $a_{i_{1}\dots i_{p+1}}\neq 0\}$. Then $d = deg\,\Phi_{\eta}$, $a_{m_{1}\dots m_{p}n}$,
elements $(k_{1},\dots, k_{p+1})\in E_{\eta}'$, the corresponding
coefficients $a_{k_{1}\dots k_{p+1}}$ and the coefficients of the
terms of total degree $d$ do not depend on the choice of the
system of $\Delta$-$\sigma$-generators $\eta$. Furthermore,
$\D\frac{a_{m_{1}\dots m_{p}n}}{2^{n}} = \Delta$-$\sigma$-$tr.deg_{K}L$.
\end{theorem}

\begin{definition}
The polynomial $\Phi_{\eta}(t_{1},\dots, t_{p+1})$ whose existence is established by Theorem 3.1
is called the difference-differential {\em (or $\Delta$-$\sigma$-)} dimension polynomial associated with
the partition  (3.1) of the basic set of derivations.
\end{definition}

The $\Delta$-$\sigma$-dimension polynomial associated with partition (3.1) has the following interpretation
as the strength of a system of difference-differential equations.

Let us consider a system of partial difference-differential equations
\begin{equation}
A_{i}(f_{1},\dots, f_{s}) = 0\hspace{0.3in}(i=1,\dots, q)
\end{equation}
over a field of of functions in $m$ real variables $x_{1},\dots, x_{m}$
($f_{1},\dots, f_{s}$ are unknown functions of $x_{1},\dots, x_{m}$).
Suppose that $\Delta = \{\delta_{1},\dots, \delta_{m}\}$ where $\delta_{i}$
is the partial differentiation $\partial/\partial x_{i}$ ($i=1,\dots, m$) and
the basic set of automorphisms $\sigma = \{\alpha_{1},\dots, \alpha_{m}\}$ consists of $m$ shifts of arguments,
$f(x_{1},\dots, x_{m})\mapsto f(x_{1},\dots, x_{i-1}, x_{i}+h_{i}, x_{i+1},\dots, x_{m})$
($1\leq i\leq m$, $h_{1},\dots, h_{m}$ are some real numbers).
Thus, we assume that the the left-hand sides of the equations in (3.2) contain
unknown functions $f_{i}$, their partial derivatives, their images under
the shifts $\alpha_{j}$ and various compositions of such shifts and partial derivations.
Furthermore, we suppose that system (3.2) is algebraic, that is, all $A_{i}(y_{1},\dots, y_{s})$
are elements of a ring of $\Delta$-$\sigma$-polynomials $K\{y_{1},\dots, y_{s}\}$ with coefficients in some
functional $\Delta$-$\sigma$-field $K$.

Let us consider a grid with equal cells of dimension
$h_{1}\times\dots\times h_{m}$ that fills the whole space ${\bf R}^{m}$.
Let us fix some node $\mathcal{P}$ say that {\em a node
$\mathcal{Q}$ has order $i$} (with respect to $\mathcal{P}$) if the
shortest path from $\mathcal{P}$ to $\mathcal{Q}$ along the edges of
the grid consists of $i$ steps (by a step we mean a path from a node
of the grid to a neighbor node along the edge between these two
nodes). Say, the orders of the nodes in the two-dimensional case are
as follows (a number near a node shows the order of this node).
\setlength{\unitlength}{1cm}
\begin{picture}(9,4.9)
\put(1.6,0.3){\line(1,0){0.7}} \put (2.3,0.30){\circle*{0.1}}
\put(2.3,0.3){\line(1,0){1}} \put (3.3,0.30){\circle*{0.1}} \put
(4.3,0.30){\circle*{0.1}} \put (6.3,0.30){\circle*{0.1}} \put
(5.3,0.30){\circle*{0.1}} \put(3.3,0.3){\line(1,0){1}}
\put(4.3,0.3){\line(1,0){1}} \put(5.3,0.3){\line(1,0){1}}
\put(6.3,0.3){\line(1,0){1}} \put (7.3,0.30){\circle*{0.1}}
\put(7.3,0.3){\line(1,0){1}} \put (8.3,0.30){\circle*{0.1}}
\put(8.3,0.3){\line(1,0){0.7}} \put(2.3,0.3){\line(0,-1){0.2}}
\put(3.3,0.3){\line(0,-1){0.2}} \put(4.3,0.3){\line(0,-1){0.2}}
\put(5.3,0.3){\line(0,-1){0.2}} \put(6.3,0.3){\line(0,-1){0.2}}
\put(7.3,0.3){\line(0,-1){0.2}} \put(8.3,0.3){\line(0,-1){0.2}}
\put(2.3,0.3){\line(0,1){1}} \put(3.3,0.3){\line(0,1){1}}
\put(4.3,0.3){\line(0,1){1}} \put(5.3,0.3){\line(0,1){1}}
\put(6.3,0.3){\line(0,1){1}} \put(7.3,0.3){\line(0,1){1}}
\put(8.3,0.3){\line(0,1){1}} \put(7.3,1.30){\circle*{0.1}}
\put(2.3,1.30){\circle*{0.1}} \put(3.3,1.30){\circle*{0.1}}
\put(4.3,1.30){\circle*{0.1}} \put(5.3,1.30){\circle*{0.1}}
\put(6.3,1.30){\circle*{0.1}} \put(8.3,1.30){\circle*{0.1}}
\put(1.6,1.3){\line(1,0){0.7}} \put(2.3,1.3){\line(1,0){1}}
\put(3.3,1.3){\line(1,0){1}} \put(4.3,1.3){\line(1,0){1}}
\put(5.3,1.3){\line(1,0){1}} \put(6.3,1.3){\line(1,0){1}}
\put(7.3,1.3){\line(1,0){1}} \put(8.3,1.3){\line(1,0){0.7}}

\put(2.3,2.3){\circle*{0.1}} \put(3.3,2.3){\circle*{0.1}}
\put(4.3,2.3){\circle*{0.1}} \put(5.3,2.3){\circle*{0.1}}
\put(6.3,2.3){\circle*{0.1}} \put(7.3,2.3){\circle*{0.1}}
\put(8.3,2.3){\circle*{0.1}}

\put(2.3,3.3){\circle*{0.1}} \put(3.3,3.3){\circle*{0.1}}
\put(4.3,3.3){\circle*{0.1}} \put(5.3,3.3){\circle*{0.1}}
\put(6.3,3.3){\circle*{0.1}} \put(7.3,3.3){\circle*{0.1}}
\put(8.3,3.3){\circle*{0.1}}

\put(2.3,4.3){\circle*{0.1}} \put(3.3,4.3){\circle*{0.1}}
\put(4.3,4.3){\circle*{0.1}} \put(5.3,4.3){\circle*{0.1}}
\put(6.3,4.3){\circle*{0.1}} \put(7.3,4.3){\circle*{0.1}}
\put(8.3,4.3){\circle*{0.1}}

\put(1.6,2.3){\line(1,0){0.7}} \put(2.3,2.3){\line(1,0){1}}
\put(3.3,2.3){\line(1,0){1}} \put(4.3,2.3){\line(1,0){1}}
\put(5.3,2.3){\line(1,0){1}} \put(6.3,2.3){\line(1,0){1}}
\put(7.3,2.3){\line(1,0){1}} \put(8.3,2.3){\line(1,0){0.7}}

\put(1.6,3.3){\line(1,0){0.7}} \put(2.3,3.3){\line(1,0){1}}
\put(3.3,3.3){\line(1,0){1}} \put(4.3,3.3){\line(1,0){1}}
\put(5.3,3.3){\line(1,0){1}} \put(6.3,3.3){\line(1,0){1}}
\put(7.3,3.3){\line(1,0){1}} \put(8.3,3.3){\line(1,0){0.7}}

\put(1.6,4.3){\line(1,0){0.7}} \put(2.3,4.3){\line(1,0){1}}
\put(3.3,4.3){\line(1,0){1}} \put(4.3,4.3){\line(1,0){1}}
\put(5.3,4.3){\line(1,0){1}} \put(6.3,4.3){\line(1,0){1}}
\put(7.3,4.3){\line(1,0){1}} \put(8.3,4.3){\line(1,0){0.7}}

\put(2.3,1.3){\line(0,1){1}} \put(3.3,1.3){\line(0,1){1}}
\put(4.3,1.3){\line(0,1){1}} \put(5.3,1.3){\line(0,1){1}}
\put(6.3,1.3){\line(0,1){1}} \put(7.3,1.3){\line(0,1){1}}
\put(8.3,1.3){\line(0,1){1}}

\put(2.3,2.3){\line(0,1){1}} \put(3.3,2.3){\line(0,1){1}}
\put(4.3,2.3){\line(0,1){1}} \put(5.3,2.3){\line(0,1){1}}
\put(6.3,2.3){\line(0,1){1}} \put(7.3,2.3){\line(0,1){1}}
\put(8.3,2.3){\line(0,1){1}}

\put(2.3,3.3){\line(0,1){1}} \put(3.3,3.3){\line(0,1){1}}
\put(4.3,3.3){\line(0,1){1}} \put(5.3,3.3){\line(0,1){1}}
\put(6.3,3.3){\line(0,1){1}} \put(7.3,3.3){\line(0,1){1}}
\put(8.3,3.3){\line(0,1){1}}

\put(2.3,4.3){\line(0,1){0.5}} \put(3.3,4.3){\line(0,1){0.5}}
\put(4.3,4.3){\line(0,1){0.5}} \put(5.3,4.3){\line(0,1){0.5}}
\put(6.3,4.3){\line(0,1){0.5}} \put(7.3,4.3){\line(0,1){0.5}}
\put(8.3,4.3){\line(0,1){0.5}}

\put(5.15,2.1){\makebox(0,0){$\mathcal{P}$}}
\put(4.15,2.1){\makebox(0,0){$1$}}
\put(6.15,2.1){\makebox(0,0){$1$}}
\put(3.15,2.1){\makebox(0,0){$2$}}
\put(7.15,2.1){\makebox(0,0){$2$}}
\put(2.15,2.1){\makebox(0,0){$3$}}
\put(8.15,2.1){\makebox(0,0){$3$}}

\put(5.15,1.1){\makebox(0,0){$1$}}
\put(5.15,3.1){\makebox(0,0){$1$}}
\put(4.15,1.1){\makebox(0,0){$2$}}
\put(6.15,1.1){\makebox(0,0){$2$}}
\put(3.15,1.1){\makebox(0,0){$3$}}
\put(7.15,1.1){\makebox(0,0){$3$}}
\put(2.15,1.1){\makebox(0,0){$4$}}
\put(8.15,1.1){\makebox(0,0){$4$}}

\put(5.15,4.1){\makebox(0,0){$2$}}
\put(5.15,0.1){\makebox(0,0){$2$}}
\put(4.15,0.1){\makebox(0,0){$3$}}
\put(6.15,4.1){\makebox(0,0){$3$}}
\put(3.15,0.1){\makebox(0,0){$4$}}
\put(7.15,4.1){\makebox(0,0){$4$}}
\put(2.15,0.1){\makebox(0,0){$5$}}
\put(8.15,4.1){\makebox(0,0){$5$}}

\put(6.15,0.1){\makebox(0,0){$3$}}
\put(7.15,0.1){\makebox(0,0){$4$}}
\put(8.15,0.1){\makebox(0,0){$5$}}

\put(6.15,4.1){\makebox(0,0){$3$}}
\put(7.15,4.1){\makebox(0,0){$4$}}
\put(8.15,4.1){\makebox(0,0){$5$}}

\put(6.15,3.1){\makebox(0,0){$2$}}
\put(7.15,3.1){\makebox(0,0){$3$}}
\put(8.15,3.1){\makebox(0,0){$4$}}

\put(4.15,4.1){\makebox(0,0){$3$}}
\put(3.15,4.1){\makebox(0,0){$4$}}
\put(2.15,4.1){\makebox(0,0){$5$}}

\put(4.15,3.1){\makebox(0,0){$2$}}
\put(3.15,3.1){\makebox(0,0){$3$}}
\put(2.15,3.1){\makebox(0,0){$4$}}
\end{picture}\\
Let us fix partition (3.1) of the set of basic derivations $\Delta$
(such a partition can be, for example,  a natural separation of (all or some)
derivations with respect to coordinates and the derivation with respect to time).
Let us consider the values of the unknown functions
$f_{1},\dots, f_{s}$ and their partial derivatives whose $i$th order
(that is the order with respect to the derivations of the set $\Delta_{i}$)
does not exceed $r_{i}$ ($1\leq i\leq p$) at the nodes whose order does not exceed $r_{p+1}$\,
($r_{1},\dots, r_{p+1}\in {\bf N}$). If $f_{1},\dots, f_{s}$ should not satisfy any
system of equations (or any other condition), these values can be chosen arbitrarily.
Because of the system (and equations obtained from the equations of the
system by partial differentiations and transformations of the form $f_{j}(x_{1},\dots,
x_{m})\mapsto f_{j}(x_{1}+k_{1}h_{1},\dots, x_{m}+k_{m}h_{m})$ with
$k_{1},\dots, k_{m}\in {\bf Z}$, $1\leq j\leq s$), the number of
independent values of the functions $f_{1},\dots, f_{s}$ and their partial
derivatives whose $i$th order does not exceed $r_{i}$ ($1\leq i\leq p$) at the
nodes of order $\leq r_{p+1}$ decreases. This number, which is a function
of $p+1$ variables $r_{1},\dots, r_{p+1}$, is the ''measure of strength'' of the system in
finite differences in the sense of A. Einstein). We denote it by $S_{r_{1},\dots, r_{p+1}}$.

Suppose that the $\Delta$-$\sigma$-ideal $J$ generated in $K\{y_{1},\dots, y_{s}\}$
by the $\Delta$-$\sigma$-polynomials $A_{1},\dots, A_{q}$ is prime (e. g., the polynomials are linear).
The we say that the system of difference-differential equations (3.2) is prime. In this case, the field
of fractions $L$ of the $\Delta$-$\sigma$-integral domain $K\{y_{1},\dots, y_{s}\}/J$ has a natural structure
of a $\Delta$-$\sigma$-field extension of $K$ generated by the finite set $\eta = \{\eta_{1},\dots, \eta_{s}\}$
where $\eta_{i}$ is the canonical image of $y_{i}$ in $K\{y_{1},\dots, y_{s}\}/J$ ($1\leq i\leq s$).
It is easy to see that the $\Delta$-$\sigma$-dimension polynomial $\Phi_{\eta}(t_{1},\dots, t_{p+1})$ of
the extension $L/K$ associated with the system of $\Delta$-$\sigma$-generators $\eta$ has the property that
$\Phi_{\eta}(r_{1},\dots, r_{p+1}) = S_{r_{1},\dots, r_{p+1}}$ for all sufficiently large
$(r_{1},\dots, r_{p+q})\in {\bf N}^{p+1}$, so this dimension polynomial is the measure of strength of
the system of difference-differential equations (3.2) in the sense of A. Einstein.

\section{Numerical polynomials of subsets of ${\bf N}^{m}\times {\bf Z}^{n}$}

\setcounter{equation}{0}

\begin{definition}
A polynomial $f(t_{1}, \dots,t_{p})$ in $p$ variables $t_{1},\dots, t_{p}$
($p\in{\bf N}, p\geq 1$) with rational coefficients is called
{\em numerical\/} if $f(t_{1},\dots, t_{p})\in {\bf Z}$
for all sufficiently large $(t_{1}, \dots, t_{p})\in{\bf Z}^{p}$. (Recall that it means that  there
exist integers $s_{1},\dots, s_{p}$ such that
$f(r_{1},\dots, r_{p})\in {\bf Z}$ as soon as
$(r_{1},\dots, r_{p})\in {\bf Z}^{p}$ and $r_{i}\geq s_{i}$ for all
$i = 1,\dots, p$.)
\end{definition}

It is clear that every polynomial with integer
coefficients is numerical.  As an example of a numerical polynomial in
$p$ variables with noninteger coefficients ($p\in {\bf N}, p\geq 1$)
one can consider a polynomial
$\D\prod_{i=1}^{p}{t_{i}\choose m_{i}}$ \, where
$m_{1},\dots, m_{p}\in{\bf N}$. (As usual, $\D{t\choose k}$ ($k\in {\bf Z},
k\geq 1$) denotes the polynomial $\D\frac{t(t-1)\dots (t-k+1)}{k!}$ in one
variable $t$, $\D{t\choose0} = 1$, and $\D{t\choose k} = 0$ if $k$ is a
negative integer.)

The following theorem proved in \cite{KLMP} gives the ''canonical''
representation of a numerical polynomial in several variables.

\begin{theorem}

Let $f(t_{1},\dots, t_{p})$ be a numerical polynomial in $p$ variables
$t_{1},\dots,  t_{p}$, and let $deg_{t_{i}}\, f = m_{i}$
($m_{1},\dots, m_{p}\in{\bf N}$). Then the
polynomial  $f(t_{1},\dots, t_{p})$ can be represented in the form

\begin{equation}
f(t_{1},\dots t_{p}) =\D\sum_{i_{1}=0}^{m_{1}}\dots \D\sum_{i_{p}=0}^{m_{p}}
{a_{i_{1}\dots i_{p}}}{t_{1}+i_{1}\choose i_{1}}\dots{t_{p}+i_{p}
\choose i_{p}}
\end{equation}

\noindent with integer coefficients $a_{i_{1}\dots i_{p}}$
($0\leq i_{k}\leq m_{k}$ for $k = 1,\dots, p$)
that are uniquely defined by the numerical polynomial.

\end{theorem}

In what follows (until the end of the section), we deal with subsets of
the set ${\bf N}^{m}\times {\bf Z}^{n}$ ($m$ and $n$ are positive integers).
Furthermore, we fix a partition of the set ${\bf N}_{m} = \{1,\dots , m\}$ into $p$ disjoint subsets ($p\geq 1$):
\begin{equation}
{\bf N}_{m} = \{1,\dots , m_{1}\}\bigcup \{m_{1}+1,\dots , m_{1}+m_{2}\}\bigcup\dots\bigcup \{m_{1}+\dots +m_{p-1}+1,\dots , m\}
\end{equation}
($m_{1}+\dots +m_{p} = m$).

If $a = (a_{1},\dots, a_{m+n})\in {\bf N}^{m}\times {\bf Z}^{n}$ we denote the numbers $\D\sum_{i=1}^{m_{1}}a_{i}$,
$\D\sum_{i=m_{1}+1}^{m_{1}+m_{2}}a_{i},\dots,$

\noindent$\D\sum_{i=m_{1}+\dots + m_{p-1} + 1}^{m}a_{i}$,
$\D\sum_{i=m+1}^{m+n}|a_{i}|$ by $ord_{1}a,\dots, ord_{p+1}a$, respectively.

\medskip

As in \cite[Section 2.5]{KLMP}, let us consider the set ${\bf Z}^{n}$ as a union
\begin{equation}
{\bf Z}^{n} = \bigcup_{1\leq j\leq 2^{n}}{\bf Z}_{j}^{(n)}
\end{equation}
where
${\bf Z}_{1}^{(n)}, \dots, {\bf Z}_{2^{n}}^{(n)}$ are all different
Cartesian products of $n$ sets each of which is either
${\bf N}$ or ${\bf Z_{-}}=\{k\in{\bf Z}\,|\,k\leq 0\}$. We assume that
${\bf Z}_{1}^{(n)} = {\bf N}^{n}$ and call ${\bf Z}_{j}^{(n)}$ the
{\em $j$th orthant} of the set ${\bf Z}^{n}$ ($1\leq j\leq 2^{n}$).
Furthermore, we consider ${\bf N}^{m}\times {\bf Z}^{n}$ as a partially
ordered set with the order $\unlhd$ such that
$(e_{1},\dots, e_{m}, f_{1},\dots, f_{n})\unlhd
(e'_{1},\dots, e'_{m}, f'_{1},\dots, f'_{n})$ if and only if
$(f_{1},\dots, f_{n})$ and $(f'_{1},\dots, f'_{n})$ belong to the same
orthant ${\bf Z}_{k}^{(n)}$ ($1\leq k\leq 2^{n}$) and the
$(m+n)$-tuple $(e_{1},\dots, e_{m}, |f_{1}|,\dots, |f_{n}|)$ is less
than $(e'_{1},\dots, e'_{m}, |f'_{1}|,\dots, |f'_{n}|)$ with respect
to the product order on ${\bf N}^{m+n}$.

In what follows, for any set $A\subseteq {\bf N}^{m}\times {\bf Z}^{n}$,
$W_{A}$ will denote the set of all elements of ${\bf N}^{m}\times {\bf Z}^{n}$
that do not exceed any element of $A$ with respect to the order
$\unlhd$.  (Thus, $w\in W_{A}$ if and only if there is no element
$a\in A$ such that $a\unlhd w$.) Furthermore, for any
$r_{1},\dots r_{p+1}\in {\bf N}$, $A(r_{1},\dots r_{p+1})$ will denote the set of
all elements $x = (x_{1},\dots, x_{m}, x'_{1},\dots, x'_{n})\in A$
such that $ord_{i}x\leq r_{i}$ ($i=1,\dots, p+1$).

\medskip

The above notation can be naturally restricted to subsets of ${\bf N}^{m}$. If $E\subseteq {\bf N}^{m}$
and $s_{1},\dots, s_{p}$ are nonnegative integers,
then $E(s_{1},\dots, s_{p})$ will denote the set of all $m$-tuples $e = (e_{1},\dots, e_{m})\in E$
such that $ord_{i}(e_{1},\dots, e_{m},0)\leq s_{i}$ for $i=1,\dots, p$.
Furthermore, we shall associate with a set $E\subseteq {\bf N}^{m} $ a set
$V_{E}\subseteq {\bf N}^{m}$ that consists of all $m$-tuples
$v = (v_{1},\dots , v_{m})\in {\bf N}$ that are not greater than or equal
to any $m$-tuple in $E$ with respect to the product order on
${\bf N}^{m}$. (Recall that the product order on  ${\bf N}^{m}$
is a partial order $\leq_P$ on ${\bf N}^{m}$ such
that $c =(c_{1}, \dots , c_{m})\leq_{P} c' =(c'_{1}, \dots , c'_{})$ if and only if
$c_{i}\leq c'_{i}$ for all $i=1, \dots , m$.   If $c\leq_{P} c'$ and
$c\neq c'$, we write $c<_{P} c'$ ).  Clearly, an element
$v=(v_{1}, \dots , v_{m})\in {\bf N}^{m}$ belongs to $V_{E}$ if and only if
for any element  $(e_{1},\dots , e_{m})\in E$ there exists
$i\in {\bf N}, 1\leq i\leq m$, such that $e_{i} > v_{i}$.

The following two theorems proved in \cite[Chapter 2]{KLMP} generalize the well-known Kolchin's
result on the numerical polynomials associated with subsets of ${\bf N}$
(see \cite[Chapter 0, Lemma 17]{K2}) and give the explicit formula for the
numerical polynomials in $p$ variables associated with a finite subset
of ${\bf N}^{m}$.

\begin{theorem}
Let $E$ be a subset of ${\bf N}^{m}$ where $m = m_{1} + \dots + m_{p}$
for some nonnegative integers  $m_{1},\dots, m_{p}$ ($p\geq 1$). Then there
exists a numerical polynomial $\omega_{E}(t_{1},\dots, t_{p})$ with the
following properties:

{\em (i)} \,  $\omega_{E}(r_{1},\dots, r_{p}) =
Card\,V_{E}(r_{1},\dots, r_{p})$
for all sufficiently large $(r_{1},\dots, r_{p})\in {\bf N}^{p}$

(as usual, $Card \, M$ denotes the number of elements of a finite set $M$).

{\em (ii)} \, The total degree of the polynomial  $\omega_{E}$ does not
exceed $m$ and $deg_{t_{i}}\omega_{E}\leq m_{i}$ for all
$i = 1,\dots, p$.

{\em (iii)} \, $deg\,\omega_{E} = m$ if and only if the $E=\emptyset$.
Then  $\omega_{E}(t_{1},\dots, t_{p}) =
\D\prod_{i=1}^{p}{t_{i}+m_{i}\choose m_{i}}$.

{\em (iv)} \, $\omega_{E}$ is a zero polynomial if and only if
$(0,\dots, 0)\in E$.

\end{theorem}

\begin{definition}
The polynomial $\omega_{E}(t_{1},\dots, t_{p})$ whose existence is
stated by Theorem 4.3 is called the dimension polynomial of the set
$E\subseteq {\bf N}^{m}$ associated with the partition
$(m_{1},\dots, m_{p})$ of $m$. If $p=1$, the polynomial $\omega_{E}$
is called the Kolchin polynomial of the set $E$.
\end{definition}

\begin{theorem}
Let $E = \{e_{1}, \dots, e_{q}\}$ be a finite subset of ${\bf N}^{m}$
where $q$ is a positive integer and $m = m_{1} + \dots + m_{p}$ for
some nonnegative integers  $m_{1},\dots, m_{p}$ ($p\geq 1$).
Let $e_{i} = (e_{i1}, \dots, e_{im})$ \, ($1\leq i\leq q$)
and for any $l\in {\bf N}$, $0\leq l\leq q$, let $\Gamma (l,q)$
denote the set of all $l$-element subsets of the set
${\bf N}_{q} = \{1,\dots, q\}$. Furthermore, for any
$\sigma \in \Gamma (l,q)$, let $\bar{e}_{\sigma j} = \max \{e_{ij} |
i\in \sigma\}$ ($1\leq j\leq m$) and
$b_{\sigma j} =\D\sum_{h\in \sigma j}\bar{e}_{\sigma h}$.
Then \begin{equation} \omega_{E}(t_{1},\dots, t_{p}) =
\D\sum_{l=0}^{q}(-1)^{l}\D\sum_{\sigma \in \Gamma (l,q)}
\D\prod_{j=1}^{p}{t_{j}+m_{j} - b_{\sigma j}\choose m_{j}}
\end{equation}
\end{theorem}

{\bf Remark.} \, It is clear that if $E$ is any subset of ${\bf N}^{m}$ and
$E^{\ast}$ is the set of all minimal elements of the set $E$ with
respect to the product order on ${\bf N}^{m}$, then the set $E^{\ast}$ is
finite and $\omega_{E}(t_{1}, \dots, t_{p}) =
\omega_{E^{\ast}}(t_{1}, \dots, t_{p})$.
Thus, Theorem 4.3 gives an algorithm that allows one to find a
numerical polynomial associated with any subset of  ${\bf N}^{m}$  (and with
a given partition of the set $\{1,\dots, m\}$): one should
first find the set of all minimal points of the subset and then apply
Theorem 4.3.

\bigskip

The following result can be obtained precisely in the same way as Theorem 3.4 of  \cite{Levin2}
(the only difference is that the proof in the mentioned paper uses Theorem 3.2 of \cite{Levin2}
in the case  $p=2$, while the proof of the theorem below should refer to the
Theorem 3.2 of \cite{Levin2} where $p$ is any positive integer).

\begin{theorem}
Let $A$ be a subset of ${\bf N}^{m}\times{\bf Z}^{n}$ and
let partition {\em (4.2)} of the set ${\bf N}_{m}$ be fixed. Then there
exists a numerical polynomial $\phi_{A}(t_{1},\dots, t_{p+1})$ in $p+1$ variables
$t_{1},\dots, t_{p+1}$ with the following properties.

{\em (i)}\, $\phi_{A}(r_{1},\dots, r_{p+1}) = Card\,W_{A}(r_{1},\dots, r_{p+1})$ for all sufficiently large

\noindent$(r_{1},\dots, r_{p+1})\in {\bf N}^{p+1}$.

\medskip

{\em (ii)}\, $deg_{t_{i}}\phi_{A}\leq m_{i}$ for $i=1,\dots, p$ and $deg_{t_{p+1}}\phi_{A}\leq n$.

\medskip

{\em (iii)}\, Let us consider a mapping $\rho: {\bf N}^{m}\times {\bf Z}^{n}\longrightarrow{\bf N}^{m+2n}$
such that\, 

\smallskip

\noindent$\rho((a_{1},\dots, a_{m+n}) =(a_{1},\dots, a_{m}, \max \{a_{m+1}, 0\}, \dots, \max \{-a_{m+1}, 0\}, $ 

\noindent$\max \{a_{m+n}, 0\},
\dots, \max \{-a_{m+n}, 0\})$.

\medskip

Let $B = \rho(A)\bigcup \{e_{1},\dots, e_{n}\}$ where $e_{i}$
($1\leq i\leq n$) is a $(m+2n)$-tuple in ${\bf N}^{m+2n}$ whose
$(m+i)$th and $(m+n+i)$th coordinates are equal to 1 and all other
coordinates are equal to 0.
Then $\phi_{A}(t_{1}, \dots, t_{p+1}) = \omega_{B}(t_{1}, \dots, t_{p+1})$ where
$\omega_{B}(t_{1},\dots,  t_{p+1})$ is the dimension polynomial of the set $B$
(see {\em Definition 4.4}) associated with the partition ${\bf N}_{m+2n} =
\{1,\dots , m_{1}\}\bigcup \{m_{1}+1,\dots , m_{1}+m_{2}\}\bigcup\dots
\bigcup \{m_{1}+\dots +m_{p-1}+1,\dots , m\}\bigcup \{m+1,\dots, m+2n\}$ of the
set ${\bf N}_{m+2n}$.

\medskip

{\em (iv)}\, If $A = \emptyset$, then
$$\phi_{A}(t_{1}, \dots, t_{p+1}) ={{t_{1}+m_{1}}\choose m_{1}}\dots {{t_{p}+m_{p}}\choose m_{p}}
\sum_{i=0}^{n}(-1)^{n-i}2^{i}{n\choose i}{{t_{p+1}+i}\choose i}.$$

\medskip

{\em (v)}\, $\phi_{A}(t_{1},\dots, t_{p+1}) = 0$ if and only if $(0, \dots, 0)\in A$.
\end{theorem}

\section{Proof of the main theorem and computation of difference-differential dimension
polynomials via characteristic sets}

\setcounter{equation}{0}

In this and next sections we prove Theorem 3.1 and give two methods of computation of difference-differential
dimension polynomials of $\Delta$-$\sigma$-field extensions. One of them is based on
constructing a characteristic set of the defining prime $\Delta$-$\sigma$-ideal of the
extension; the other approach is the computation of the dimension polynomial of the
associated module of K\"aller differentials via (generalized) Gr\"obner basis method.
Both approaches use total term orderings with respect to several orders defined by
partitions of the basic sets of derivations and automorphisms.

In what follows we use the notation and conventions introduced at the beginning of section 3.
In particular, we assume that partition (3.1) of the set of basic derivations
$\Delta = \{\delta_{1},\dots, \delta_{m}\}$ is fixed.

\medskip
Let us consider $p+1$ total orderings $<_{1}, \dots, <_{p}, <_{\sigma}$ of the
set of power products $\Lambda$ such that

\medskip

\noindent$\lambda  = \delta_{1}^{k_{1}}\dots \delta_{m}^{k_{m}}
\alpha_{1}^{l_{1}}\dots \alpha_{n}^{l_{n}}
<_{i} \lambda'  = \delta_{1}^{k'_{1}}\dots
\delta_{m}^{k'_{m}}\alpha_{1}^{l'_{1}}\dots \alpha_{n}^{l'_{n}}$  ($1\leq i\leq p$)  if and only if

\medskip

\noindent $(ord_{i}\lambda, ord\,\lambda, ord_{1}\lambda,\dots, ord_{i-1}\lambda, ord_{i+1}\lambda, \dots, ord_{p}\lambda,
ord_{\sigma}\lambda,  k_{m_{1}+\dots + m_{i-1}+1},\dots,$

 \medskip

\noindent$k_{m_{1}+\dots +m_{i}},\, k_{1},\dots, k_{m_{1}+\dots +
m_{i-1}}, k_{m_{1}+\dots + m_{i}+1},\dots, k_{m}, |l_{1}|,\dots, |l_{n}|, l_{1},\dots, l_{n})$ is

\medskip

\noindent less than
$(ord_{i}\lambda', ord\,\lambda', ord_{1}\lambda',\dots, ord_{i-1}\lambda', ord_{i+1}\lambda', \dots, ord_{p}\lambda', ord_{\sigma}\lambda',$

\medskip

\noindent$ k'_{m_{1}+\dots + m_{i-1}+1},\dots, k'_{m_{1}+\dots +m_{i}},\, k'_{1},\dots, k'_{m_{1}+\dots +
m_{i-1}}, k'_{m_{1}+\dots + m_{i}+1},\dots,$

\medskip

\noindent$k'_{m}, |l'_{1}|,\dots, |l'_{n}|, l'_{1},\dots, l'_{n})$ with respect to the lexicographic order on ${\bf N}^{m+2n+p+2}$.

\medskip

Similarly, $\lambda <_{\sigma} \lambda'$ if and only if $(ord_{\sigma}\lambda, ord\,\lambda, ord_{1}\lambda,\dots, ord_{p}\lambda,
|l_{1}|,\dots, |l_{n}|,$

\medskip

\noindent$l_{1},\dots, l_{n}, k_{1},\dots, k_{m})$ is less than the corresponding $(m+2n+p+2)$-tuple for $\lambda'$
with respect to the lexicographic order on ${\bf N}^{m+2n+p+2}$.

\medskip

Two elements $\lambda_{1} = \delta_{1}^{k_{1}}\dots \delta_{m}^{k_{m}}
\alpha_{1}^{l_{1}}\dots \alpha_{n}^{l_{n}}$ and
$\lambda_{2} = \delta_{1}^{r_{1}}\dots \delta_{m}^{r_{m}}\alpha_{1}^{s_{1}}
\dots \alpha_{n}^{s_{n}}$ in $\Lambda$ are called {\em similar\/}, if
the $n$-tuples $(l_{1}, \dots, l_{n})$ and $(s_{1}, \dots, s_{n})$
belong to the same orthant of ${\bf Z}^{n}$ (see (4.3)\,). In this case we
write $\lambda_{1}\sim \lambda_{2}$. We say that $\lambda_{1}$ {\em divides} $\lambda_{2}$
(or $\lambda_{2}$ is a {\em multiple} of $\lambda_{1}$) and write $\lambda_{1}|\lambda_{2}$
if $\lambda_{1}\sim \lambda_{2}$ and there exists
$\lambda \in \Lambda$ such that $\lambda\sim \lambda_{1}$ and $\lambda_{2} = \lambda\lambda_{1}$.

\medskip

Let $K$ be a difference-differential field ($Char\,K = 0$) with the basic set $\Delta\bigcup\sigma$ described above
and let partition (3.1) of the set $\Delta$ be fixed. Let $K\{y_{1},\dots, y_{s}\}$ be the ring of
$\Delta$-$\sigma$-polynomials over $K$ and let $\Lambda Y$ denote
the set of all elements $\lambda y_{i}$ ($\lambda\in \Lambda$, $1\leq i\leq s$) called {\em terms}.
Note that as a ring, $K\{y_{1},\dots, y_{s}\} = K[\Lambda Y]$. Two terms $u=\lambda y_{i}$ and $v=\lambda' y_{j}$
are called {\em similar} if $\lambda$ and $\lambda'$ are similar; in this case we write $u\sim v$.
If $u = \lambda y_{i}$ is a term and $\lambda'\in \Lambda$, we say that $u$ is similar to $\lambda'$ and write
$u\sim \lambda'$ if $\lambda\sim \lambda'$. Furthermore, if $u, v\in \Lambda Y$, we say that $u$ {\em divides} $v$ or
{\em $v$ is a multiple of $u$}, if $u=\lambda y_{i}$, $v=\lambda' y_{i}$ for some $y_{i}$ and $\lambda|\lambda'$.

\medskip

Let us consider $p+1$ orders $<_{1},\dots, <_{p}, <_{\sigma}$ on the set $\Lambda Y$ that
correspond to the orders on the semigroup $\Lambda$ (we use the same symbols for the orders on $\Lambda$
and $\Lambda Y$). These orders  are defined as follows: $\lambda y_{j} <_{i}$ (or $<_{\sigma}$) $\lambda' y_{k}$ if and only if
$\lambda <_{i}$ (respectively, $<_{\sigma}$)$\lambda'$ in $\Lambda$ or $\lambda = \lambda'$ and $j < k$ ($1\leq i\leq p,\, 1\leq j, k\leq s$).

\bigskip

The order of a term $u = \lambda y_{k}$ and its orders with
respect to the sets $\Delta_{i}$ ($1\leq i\leq p$) and $\sigma$
are defined as the corresponding orders of $\lambda$ (we use the same notation
$ord\,u$, $ord_{i}u$, and $ord_{\sigma}$ for the corresponding orders).

\medskip

If $A\in K\{y_{1},\dots, y_{s}\}\setminus K$ and $1\leq k\leq
p$, then the highest with respect to $<_{k}$ term
that appears in $A$ is called the {\em $k$-leader} of the
$\Delta$-$\sigma$-polynomial $A$. It is denoted by $u_{A}^{(k)}$. The highest
term of $A$ with respect to $<_{\sigma}$ is called the {\em $\sigma$-leader} of
$A$; it is denoted by $v_{A}$.

\smallskip

If $A$ is written as a polynomial in  $v_{A}$,
$A = I_{d}{(v_{A})}^{d} + I_{d-1}{(v_{A})}^{d-1} + \dots +
I_{0}$, where all terms of $I_{0},\dots, I_{d}$ are less than $v_{A}$ with respect to $<_{\sigma}$,
then $I_{d}$ is said to be the {\em initial\/} of $A$. The partial derivative of $A$ with respect to $v_{A}$,
$\partial A/\partial v_{A} = dI_{d}(v_{A})^{d-1} + (d-1)I_{d-1}{(v_{A})}^{d-2} + \dots + I_{1}$, is called the
{\em separant\/} of $A$. The leading coefficient and the separant of a
$\Delta$-$\sigma$-polynomial $A$ are denoted by $I_{A}$ and $S_{A}$, respectively.

\smallskip

If $A, B\in K\{y_{1},\dots, y_{s}\}$, then $A$ is said to have
lower rank than $B$ (we write $rk\,A < rk\,B$) if either $A\in K$, $B\notin K$, or
$(v_{A}, deg_{v_{A}}A, ord_{1}u_{A}^{(1)},\dots,
ord_{p}u_{A}^{(p)})$ is less than
$(v_{B}, deg_{v_{B}}B, ord_{1}u_{B}^{(1)},\dots, ord_{p}u_{B}^{(p)})$ with respect to the
lexicographic order ($v_{A}$ and $v_{B}$ are compared with respect to $<_{\sigma}$).

If the vectors are equal (or $A, B\in K$) we say that $A$
and $B$ are of the same rank and write $rk\,A = rk\, B$.

\begin{definition}
If $A, B\in K\{y_{1},\dots, y_{s}\}$, then $B$ is said to be
reduced with respect to $A$ if

(i) $B$ does not contain terms $\lambda v_{A}$ such that $\lambda\sim v_{A}$,
$\lambda_{\Delta}\neq 1$, and  $ord_{i}(\lambda u_{A}^{(i)})\leq ord_{i}u_{B}^{(i)}$
for $i=1,\dots, p$.

\medskip

(ii) If $B$ contains a term $\lambda v_{A}$, where $\lambda\sim v_{A}$,
$\lambda_{\Delta} = 1$, then either there exists
$j, \,  1\leq j\leq p$, such that  $ord_{j} u_{B}^{(j)} < ord_{j}(\lambda u_{A}^{(j)})$
or $ord_{j}(\lambda u_{A}^{(j)}) \leq ord_{j} u_{B}^{(j)}$
for all $j = 1,\dots, p$ and $deg_{\lambda v_{A}}B < deg_{v_{A}}A$.
\end{definition}

If $B\in K\{y_{1},\dots, y_{s}\}$, then $B$ is said to
be {\em reduced with respect to a set\, $\Sigma \subseteq
K\{y_{1},\dots, y_{s}\}$} if $B$ is reduced with respect to every
element of $\Sigma$.
\medskip
A set  $\Sigma \subseteq K\{y_{1},\dots, y_{s}\}$ is called
{\em autoreduced} if $\Sigma \bigcap K =
\emptyset$ and every element of $\Sigma$ is reduced with respect
to any other element of this set.

\medskip

The proof of the following lemma can be found in ~\cite[Chapter 0, Section 17]{K2}.

\begin{lemma}
Let $A$ be any infinite subset of the set  ${\bf N}^{m}\times{\bf N}_{n}$
($m,n\in {\bf N}$, $n\geq 1$). Then there exists an infinite sequence of
elements of $A$, strictly increasing relative to the product order,
in which every element has the same projection on ${\bf N}_{n}$.
\end{lemma}

This lemma immediately implies the following statement that will be used
below.

\begin{lemma}
Let $S$ be any infinite set of terms $\lambda y_j$ ($\lambda\in \Lambda,
1\leq j\leq s$) in the ring $K\{y_{1},\dots, y_{s}\}$.   Then there exists
an index $j$ ($1\leq j\leq s$) and an infinite sequence of terms
$\lambda_{1}y_{j}, \lambda_{2}y_{j}, \dots, \lambda_{k}y_{j},\dots $ such that
$\lambda_{k}|\lambda_{k+1}$ for every $k=1, 2, \dots $.
\end{lemma}

\begin{proposition}
Every autoreduced set is finite.
\end{proposition}

PROOF.\,
Suppose that $\Sigma$ is an infinite autoreduced subset
of $K\{y_{1},\dots, y_{s}\}$. Then  $\Sigma$ must contain an infinite
set  $\Sigma'\subseteq \Sigma$ such that all $\Delta$-$\sigma$-polynomials from
$\Sigma'$ have different $\sigma$-leaders similar to each other. Indeed, if it is not so, then there
exists an infinite set $\Sigma_{1}\subseteq \Sigma$ such that all
$\Delta$-$\sigma$-polynomials in $\Sigma_{1}$ have the same $\sigma$-leader $v$.
By Lemma 5.2, the infinite set
$\{(ord_{1}u_{A}^{(1)}, \dots, ord_{p}u_{A}^{(p)}) | A\in \Sigma_{1}\}$
contains a nondecreasing infinite sequence 
$(ord_{1}u_{A_{1}}^{(1)}, \dots, ord_{p}u_{A_{1}}^{(p)})
\leq_{P} (ord_{1}u_{A_{2}}^{(1)},$

\noindent$ \dots, ord_{p}u_{A_{2}}^{(p)}) \leq_{P} \dots $
($A_{1}, A_{2}, \dots \in \Sigma_{1}$ and $\leq_{P}$ denotes
the product order on ${\bf N}^{p}$) such that $v_{A_{i}}|v_{A_{i+1}}$. Since the sequence
$\{deg_{v}A_{i} | i = 1, 2, \dots \}$ cannot be strictly decreasing,
there exists two indices $i$ and $j$ such that $i < j$ and
$deg_{v_{A_{i}}}A_{i} \leq deg_{v_{A_{j}}}A_{j}$. We obtain that $A_{j}$ is reduced with
respect to  $A_{i}$ that contradicts the fact that $\Sigma$ is an
autoreduced set.

Thus, we can assume that all $\Delta$-$\sigma$-polynomials of our infinite
autoreduced set $\Sigma$ have distinct $\sigma$-leaders. By Lemma 5.3, there
exists an infinite  sequence
$B_{1}, B_{2}, \dots $ \, of elements of $\Sigma$ such that
$v_{B_{i}} | v_{B_{i+1}}$ for all $i = 1, 2, \dots $.
Let $k_{ij} = ord_{\sigma}v_{B_{i}}$ and  $l_{ij} = ord_{j}u_{B_{i}}^{(j)}$
($1\leq j\leq p$). Obviously, $l_{ij}\geq k_{ij}$ ($i = 1, 2,\dots ;
j = 1,\dots, p$), so that $\{(l_{i1}-k_{i1}, \dots, l_{ip}-k_{ip}) |
i =1, 2, \dots \}\subseteq {\bf N}^{p}$. By Lemma 5.2, there exists
an infinite sequence of indices $i_{1} < i_{2} < \dots $ such that
$(l_{i_{1}1}-k_{i_{1}1},\dots, l_{i_{1}p}-k_{i_{1}p}) \leq_{P}
(l_{i_{2}1}-k_{i_{2}1},\dots, l_{i_{2}p}-k_{i_{2}p}) \leq_{P} \dots $.
Then for any $j = 1,\dots, p$, we have
$ord_{j}\,(\, \frac{v_{B_{i_{2}}}}{v_{B_{i_{1}}}}
u_{B_{i_{1}}}^{(j)}) = k_{i_{2}j} - k_{i_{1}j} + l_{i_{1}j} \leq
k_{i_{2}j} + l_{i_{2}j} - k_{i_{2}j} =
l_{i_{2}j} = ord_{j}u_{B_{i_{2}}}^{(j)}$, so that $B_{i_{2}}$ contains a term
$\lambda v_{B_{i_{1}}} =  v_{B_{i_{2}}}$ such that $\lambda \neq 1$
and $ord_{j}(\lambda u_{B_{i_{1}}}^{(j)}) \leq ord_{j} u_{B_{i_{2}}}^{(j)}$
for $j = 1, \dots, p$. Thus, the $\Delta$-$\sigma$-polynomial $B_{i_{2}}$ is reduced
with respect to $B_{i_{1}}$ that contradicts the fact that $\Sigma$ is
an autoreduced set.

\bigskip

Throughout the rest of the paper, while considering autoreduced sets in
the ring $K\{y_{1},\dots, y_{s}\}$ we always
assume that their elements are arranged in order of increasing rank.
(Therefore, if we consider an autoreduced set of $\Delta$-polynomials
$\Sigma = \{A_{1},\dots, A_{r}\}$, then $rk\,A_{1}< \dots < rk\,A_{r}$).

\begin{theorem} Let $\Sigma = \{A_{1},\dots, A_{d}\}$ be an autoreduced set in the ring
$R = K\{y_{1},\dots, y_{s}\}$ and let $I_{k}$ and $S_{k}$ denote the initial and separant of $A_{k}$,
respectively. Furthermore, let $I(\Sigma) = \{X\in K\{y_{1},\dots, y_{s}\}\,|\,X=1$ or
$X$ is a product of finitely many elements of the form $\gamma(I_{k})$ and $\gamma'(S_{k})$
where $\gamma, \gamma'\in \Lambda_{\sigma}\}$.
Then for any $\Delta$-$\sigma$-polynomial $B$,
there exist $B_{0}\in K\{y_{1},\dots, y_{s}\}$ and $J\in I(\Sigma)$
such that $B_{0}$ is reduced with respect to $\Sigma$
and $JB\equiv B_{0} \,(mod [\Sigma])$ (that is, $JB-B_{0}\in [\Sigma]$).
\end{theorem}

PROOF.\, 
If $B$ is reduced with respect to $\Sigma$, the statement is obvious
(one can set $B_{0}=B$). Suppose that $B$ is not reduced with respect to
$\Sigma$.
Let $u^{(j)}_i$ and $v_i$ ($1\leq j\leq p,\, 1\leq i\leq d$) be the leaders of the element $A_i$ relative to
the orders $<_{j}$ and $<_{\sigma}$, respectively.
In what follows, a term $w_H$, that appears in a $\Delta$-$\sigma$-polynomial $H\in R$, will be
called a $\Sigma$-leader of $H$ if $w_H$ is the greatest (with respect to
$<_{\sigma}$) term among all terms $\lambda v_i$ ($1\leq i\leq d$) such that
$\lambda\sim v_{i}$, $\lambda v_{i}$ appears in $H$ and either $\lambda_{\Delta}\neq 1$ and
$ord_{j}(\lambda u^{(j)}_{i})\leq ord_{j}u^{(j)}_H$ for $j=1,\dots, p$, or $\lambda_{\Delta} = 1$,
$ord_{j}(\lambda u^{(j)}_{i})\leq ord_{j}u^{(j)}_H$ ($1\leq j\leq p$),
and $deg_{v_{i}}A_{i} \leq deg_{\lambda v_{i}}H$.

Let $w_B$ be the $\Sigma$-leader of $B$. Then the $\Delta$-$\sigma$-polynomial $B$ can be written as $B = B'w_{B}^{r} + B''$
where $B'$ does not contain $w_B$ and $deg_{w_{B}}B'' < r$. Let $w_{B}=\lambda v_i$
for some $i$ ($1\leq i\leq d$) and for some $\lambda \in \Lambda$, $\lambda\sim w_{B}\sim v_i$,
such that $ord_{j}(\lambda u^{(j)}_{i})\leq ord_{j}u^{(j)}_B$ for $j=1,\dots, p$.
Without loss of generality we may assume that $i$ corresponds to the maximum (with respect to the order
$<_\sigma$) $\sigma$-leader $v_i$ in the set of all $\sigma$-leaders of
elements of $\Sigma$.

Suppose, first, that $\lambda_{\Delta}\neq 1$ (and $ord_{j}(\lambda u^{(j)}_{i})\leq ord_{j}u^{(j)}_B$ for $j=1,\dots, p$).
Then $\lambda_{\Delta}A_{i} - S_{i}\lambda_{\Delta} v_{i}$ has lower rank than $\lambda_{\Delta} v_{i}$, hence
$T = \lambda A_{i} - \lambda_{\sigma}(S_{i})\lambda v_{i}$ has lower rank than $\lambda v_{i} = w_{B}$. Also,
$(\lambda_{\sigma}(S_{i}))^{r}B = (\lambda_{\sigma}(S_{i})\lambda v_{i})^{r}B' + (\lambda_{\sigma}(S_{i}))^{r}B'' =
(\lambda A_{i} - T)^{r}B' + (\lambda_{\sigma}(S_{i}))^{r}B''$. Setting $B^{(1)} = B'(-T)^{r} + (\lambda_{\sigma}(S_{i}))^{r}B''$ we
obtain that $B^{(1)}\equiv B\, mod [\Sigma]$, $B^{(1)}$ is reduced with respect to $\Sigma$ in the sense of Definition 5.1, $B^{ (1)}$
does not contain any $\Sigma$-leader, which is greater than $w_{B}$ with respect to $<_{\sigma}$, and $deg_{w_{B}} B^{(1)} < r$.

Now let $\lambda_{\Delta} = 1$, $ord_{j}(\lambda u^{(j)}_{i})\leq ord_{j}u^{(j)}_B$ ($1\leq j\leq p$),
and $r_{i} < r$ where $r_{i} = deg_{v_{i}}A_{i}$. Then the $\Delta$-$\sigma$-polynomial
$(\lambda I_{i})B - w_{B}^{r-r_{i}}(\lambda A_{i})B'$
has all the properties of $B^{(1)}$ mentioned above. Repeating the described procedure,
we arrive at a desired $\Delta$-$\sigma$-polynomial $B_{0}$
that does not contain any $\Sigma$-leader (so it is reduced with respect to $\Sigma$) and satisfies the condition
$JB\equiv B_{0} \,(mod [\Sigma])$ where $J=1$ or $J$ is a product of finitely many elements of the form $\gamma(I_{k})$
and $\gamma'(S_{k})$  ($\gamma, \gamma'\in \Lambda_{\sigma}$).

\bigskip

With the notation of the last theorem, we say that the $\Delta$-$\sigma$-polynomial $B$
{\em reduces to $B_{0}$} modulo $\Sigma$.

\begin{definition}
Let  $\Sigma = \{A_{1},\dots,A_{d}\}$ and
$\Sigma' = \{B_{1},\dots,B_{e}\}$ be two autoreduced sets in the ring of
differential polynomials  $K\{y_{1},\dots, y_{s}\}$. An autoreduced set
$\Sigma$ is said to have lower rank than $\Sigma'$ if one of the following
two cases holds:

(1) There exists $k\in {\bf N}$ such that $k\leq \min \{d,e\}$,
$rk\,A_{i}=rk\,B_{i}$ for $i=1,\dots,k-1$ and  $rk\,A_{k} < rk\,B_{k}$.

(2) $d>e$ and  $rk\,A_{i}=rk\,B_{i}$ for $i=1,\dots,e$.

If $d=e$ and $rk\,A_{i}=rk\,B_{i}$ for $i=1,\dots,d$, then $\Sigma$ is
said to have the same rank as $\Sigma'$.
\end{definition}

\begin{proposition}
In every nonempty family of autoreduced sets of differential polynomials
there exists an autoreduced set of lowest rank.
\end{proposition}

PROOF.\, Let $\Phi$ be any nonempty family of autoreduced sets
in the ring  $K\{y_{1},\dots, y_{s}\}$. Let us inductively define an infinite
descending chain of subsets of $\Phi$ as follows:
$\Phi_{0}=\Phi$, $\Phi_{1}=\{\Sigma \in \Phi_{0} | \Sigma$
contains at least one element and the first element of $\Sigma$ is of lowest
possible rank\}, \dots , $\Phi_{k}=\{\Sigma \in \Phi_{k-1} | \Sigma$ contains
at least $k$ elements and the $k$th element of $\Sigma$ is of lowest
possible rank\}, \dots . It is clear that if $A$ and $B$ are any two
$\Delta$-$\sigma$-polynomials in the same set  $\Phi_{k}$, then $v_{A} =
v_{B}$, $deg_{v_{A}}A =  deg_{v_{B}}B$, and
$ord_{i}u_{A}^{(i)} = ord_{i}u_{B}^{(i)}$ for $i = 1,\dots, p$.
Therefore, if all sets $\Phi_{k}$ are nonempty,
then the set \{$A_{k}|A_{k}$ is the $k$th element of some autoreduced set
in $\Phi_{k}$\} would be an infinite autoreduced set, and this would
contradict Lemma 5.3. Thus, there is
the smallest positive integer $k$ such that $\Phi_{k}$ is empty.
It is clear that every element of $\Phi_{k-1}$ is an autoreduced set  of
lowest rank in $\Phi$.

\bigskip

Let $J$ be any ideal of the ring   $K\{y_{1},\dots, y_{s}\}$. Since the set
of all autoreduced subsets of $J$ is not empty (if $A\in J$, then $\{A\}$
is an autoreduced subset of $J$), the last statement shows that the ideal $J$
contains an autoreduced subset of lowest rank. Such an autoreduced set is called
a {\em characteristic set} of the ideal $J$.

\begin{proposition}
Let $\Sigma = \{A_{1}, \dots , A_{d}\}$ be a characteristic set of a $\Delta$-$\sigma$-ideal
$J$ of the ring  $R = K\{y_{1},\dots, y_{s}\}$.  Then
an element $B\in R$ is reduced with respect to the set $\Sigma$ if and only
if $B = 0$.
\end{proposition}

PROOF.\,  First of all, note that if $B\neq 0$ and
$rk\,B < rk\,A_{1}$, then $rk\,\{B\} < rk\,\Sigma$
that contradicts the fact that $\Sigma$ is a characteristic set of the ideal
$J$. Let $rk\,B > rk\,A_{1}$ and let $A_{1},\dots, A_{j}$ ($1\leq j\leq d$)
be all elements of $\Sigma$ whose rank is lower that the rank of $B$.
Then the set $\Sigma' = \{A_{1},\dots, A_{j}, B\}$ is autoreduced. Indeed,
by the conditions of the theorem, $\Delta$-$\sigma$-polynomials
$A_{1},\dots, A_{j}$ are reduced with respect to each other and $B$ is
reduced with respect to the set $\{A_{1},\dots, A_{j}\}$.  Furthermore, each
$A_{i}$  ($1\leq i\leq j$) is reduced with respect to $B$ because
$rk\,A_{i} < rk\,B$. Since $rk\,\Sigma' < rk\,\Sigma$, $\Sigma$ is not a
characteristic set of $J$ that contradicts the conditions of the theorem.
Thus, $B = 0$.

\bigskip

Now we can introduce the concept of a coherent autoreduced set
of a linear $\Delta$-$\sigma$-ideal of $K\{y_{1},\dots, y_{s}\}$ (that is, a $\Delta$-$\sigma$-ideal generated by a finite set
of linear $\Delta$-$\sigma$-polynomials) in the same way as it is defined in the case of difference polynomials (see ~\cite[Section 6.5]{KLMP}):
an autoreduced set $\Sigma = \{A_{1},\dots, A_{d}\}\subseteq K\{y_{1},\dots, y_{s}\}$ consisting of linear $\Delta$-$\sigma$-polynomials is called
{\em coherent} if it satisfies the following two conditions:

(i)\, $\lambda A_{i}$ reduces to zero modulo $\Sigma$ for any $\lambda\in \Lambda, \,1\leq i\leq d$.

(ii)\, If $v_{A_{i}}\sim v_{A_{j}}$ and $w = \lambda v_{A_{i}} = \lambda'v_{A_{j}}$, where $\lambda\sim\lambda'\sim v_{A_{i}}\sim v_{A_{j}}$,
then the $\Delta$-$\sigma$-polynomial $(\lambda'I_{A_{j}})(\lambda A_{i}) - (\lambda I_{A_{i}})(\lambda'A_{j})$ reduces to zero modulo $\Sigma$.

The following two propositions can be proved precisely in the same way as the corresponding statements for difference polynomials,
see \cite[Theorem 6.5.3 and Corollary 6.5.4]{KLMP}).

\begin{proposition}
Any characteristic set of a linear $\Delta$-$\sigma$-ideal of the ring of $\Delta$-$\sigma$-polynomials $K\{y_{1},\dots, y_{s}\}$
is a coherent autoreduced set. Conversely, if $\Sigma$ is a coherent autoreduced set in $K\{y_{1},\dots, y_{s}\}$ consisting of linear
$\Delta$-$\sigma$-polynomials, then $\Sigma$ is a characteristic set of the linear $\Delta$-$\sigma$-ideal $[\Sigma]$.
\end{proposition}
\begin{proposition}
Let us consider a partial order $\preccurlyeq$ on $K\{y_{1},\dots, y_{s}\}$ such that $A\preccurlyeq B$ if and only if $v_{A}|v_{B}$. Let $A$
be a linear $\Delta$-$\sigma$-polynomial in $K\{y_{1},\dots, y_{s}\}\setminus K$. Then the set of all minimal with respect to $\preccurlyeq$
elements of the set $\{\lambda A\,|\,\lambda\in\Lambda\}$ is a characteristic set of the $\Delta$-$\sigma$-ideal $[A]$.
\end{proposition}

Now we are ready to prove Theorem 3.1.

PROOF.\,  Let $L = K\langle\eta_{1},\dots, \eta_{s}\rangle$ be a $\Delta$-$\sigma$-field extension of $K$ generated by a finite set
$\eta = \{\eta_{1},\dots, \eta_{s}\}$.  Then there exists a natural $\Delta$-$\sigma$-homomorphism $\Upsilon_{\eta}$ of the ring of $\Delta$-$\sigma$-polynomials
$K\{y_{1},\dots, y_{s}\}$ onto the $\Delta$-$\sigma$-subring $K\{\eta_{1},\dots,\eta_{s}\}$ of $L$ such that $\Upsilon_{\eta}(a) = a$ for any
$a\in K$ and $\Upsilon_{\eta}(y_{j}) = \eta_{j}$ for $j = 1,\dots, n$.
(If $A\in K\{y_{1},\dots, y_{s}\}$, then $\Upsilon_{\eta}(A)$ is called the {\em value\/} of $A$ at $\eta$ and
it is denoted by $A(\eta)$.) Obviously, the kernel $P$ of the
$\Delta$-$\sigma$-homomorphism  $\Upsilon_{\eta}$ is a prime $\Delta$-$\sigma$-ideal of the ring
$K\{y_{1},\dots, y_{s}\}$. This ideal is called the {\em defining\/} ideal
of $\eta$ over $K$ or the defining ideal of the $\Delta$-$\sigma$-field extension
$G = K\langle \eta_{1},\dots,\eta_{s}\rangle$.  It is easy to see that if
the quotient field $Q$ of the factor ring $\bar{R} =
K\{y_{1},\dots, y_{s}\}/P$ is considered as a $\Delta$-$\sigma$-field (where
$\delta(\frac{f}{g}) =
\frac{f\delta(g)-f\delta(g)}{g^2}$ and $\tau(\frac{f}{g}) = \frac{\tau(f)}{\tau(g)}$ for any $f, g\in \bar{R}$,
$\delta\in \Delta,\, \tau\in \sigma^{\ast}$), then this quotient field is naturally $\Delta$-$\sigma$-isomorphic to the field $L$.
The $\Delta$-$\sigma$-isomorphism of $Q$ onto $L$ is identity on $K$ and maps
the images of the $\Delta$-$\sigma$-indeterminates $y_{1},\dots, y_{s}$ in the
factor ring $\bar{R}$ onto the elements $\eta_{1},\dots, \eta_{s}$, respectively.

Let  $\Sigma = \{A_{1},\dots, A_{d}\}$ be a characteristic set of the defining $\Delta$-$\sigma$-ideal $P$.
For any $r_{1},\dots, r_{p+1}\in {\bf N}$, let us set $U_{r_{1}\dots r_{p+1}} =
\{u\in \Lambda Y| ord_{i}u\leq r_{i}$ for $i = 1,\dots, p+1$ and either $u$ is
not a multiple of any $v_{A_{i}}$ or for every $\lambda\in\Lambda, A\in \Sigma$ such that
$u = \lambda v_{A}$ and $\lambda\sim v_{A}$, there exists $j\in \{1,\dots, p\}$
such that $ord_{j}(\lambda u_{A}^{(j)}) > r_{j}\}$.
(For shortness, here and below we sometimes write $\leq_{i+1}$ instead of $\leq_{\sigma}$.)

We are going to show that the set $\bar{U}_{r_{1}\dots r_{p}} =
\{u(\eta)| u\in U_{r_{1}\dots r_{p}}\}$ is a transcendence basis of the
field $K(\D\bigcup_{j=1}^{n} \Lambda(r_{1},\dots, r_{p+1})\eta_{j})$ over $K$.

First of all, let us show that the set  $\bar{U}_{r_{1}\dots r_{p+1}}$ is
algebraically independent over $K$. Let $g$ be a polynomial in $k$
variables ($k\in {\bf N},  k \geq 1$) such that
$g(u_{1}(\eta),\dots, u_{k}(\eta)) = 0$ for some elements
$u_{1},\dots, u_{k}\in U_{r_{1}\dots r_{p+1}}$. Then the $\Delta$-$\sigma$-polynomial
$\bar{g} = g(u_{1},\dots, u_{k})$ is reduced with respect to
$\Sigma$. (Indeed, if $g$ contains a term $u = \lambda v_{A_{i}}$  with $\lambda\in\Lambda,\, \lambda\sim v_{A_{i}}$
($1\leq i\leq d$), then there exists $k\in \{1,\dots, p\}$ such that $ord_{k}(\lambda u_{A_{i}}^{(k)}) > r_{k}\geq ord_{k}u_{\bar{g}}^{(k)}$).
Since $\bar{g}\in P$, Proposition 5.8 implies that $\bar{g} = 0$. Thus, the set
$\bar{U}_{r_{1}\dots r_{p}}$ is algebraically independent over $K$.

Now, let us prove that every element
$\lambda \eta_{j}$ ($1\leq j\leq s, \lambda \in \Lambda(r_{1},\dots, r_{p+1})$)
is algebraic over the field $K(\bar{U}_{r_{1},\dots, r_{p+1}})$. Let
$\lambda \eta_{j} \notin \bar{U}_{r_{1},\dots, r_{p+1}}$ (if
$\lambda \eta_{j} \in \bar{U}_{r_{1},\dots, r_{p+1}}$, the statement is obvious).
Then  $\lambda y_{j} \notin U_{r_{1},\dots, r_{p}}$ whence  $\lambda y_{j}$
is equal to some term of the form $\lambda'v_{A_{i}}$ where $\lambda'\in\Lambda$, $\lambda\sim v_{A_{i}}$ ($1\leq i\leq d$), and
$ord_{k}(\lambda'u_{A_{i}}^{(k)})\leq r_{k}$ for $k = 1, \dots, p$.
Let us represent $A_{i}$ as a polynomial in $v_{A_{i}}$:
$A_{i} = I_{0}{(v_{A_{i}})}^{e} + I_{1}{(v_{A_{i}})}^{e-1} +\dots +
I_{e}$, where $I_{0}, I_{1},\dots I_{e}$ do not contain $v_{A_{i}}$
(therefore, all terms in these $\Delta$-$\sigma$-polynomials are lower than
$v_{A_{i}}$ with respect to the order $<_{\sigma}$). Since $A_{i}\in P$,
\begin{equation}
A_{i}(\eta) =  I_{0}(\eta){(v_{A_{i}}(\eta))}^{e} +
I_{1}(\eta){(v_{A_{i}}(\eta))}^{e-1} +\dots + I_{e}(\eta) = 0
\end{equation}
It is easy to see that  the $\Delta$-$\sigma$-polynomials $I_{0}$ and
$S_{A_{i}} = \partial A_{i}/\partial v_{A_{i}}$ are reduced with
respect to any element of the set $\Sigma$. Applying Proposition 5.8 we obtain
that $I_{0}\notin P$ and $S_{A_{i}}\notin P$ whence $I_{0}(\eta) \neq 0$
and $S_{A_{i}}(\eta)\neq 0$. Now, if we apply $\lambda'$ to
both sides of equation (5.1), the resulting equation will show that the
element $\lambda'v_{A_{i}}(\eta) = \lambda\eta_{j}$ is algebraic
over the field $K(\{\bar{\lambda}\eta_{l} | ord_{i}\bar{\lambda}\leq r_{i}$
for $i=1,\dots, p+1;  1\leq l\leq s$, and
$\bar{\lambda}y_{l} <_{\sigma} \lambda'u_{A_{i}}^{(1)} = \lambda y_{j}\}$.
Now, the induction on the set of terms $\Lambda Y$ ordered by the
relation $<_{\sigma}$ completes the proof of the fact that
$\bar{U}_{r_{1}\dots r_{p+1}}(\eta)$ is a transcendence basis of the
field $K(\D\bigcup_{j=1}^{s} \Lambda(r_{1},\dots, r_{p+1})\eta_{j})$ over $K$.

Let $U_{r_{1}\dots r_{p+1}}^{(1)} = \{u\in \Lambda Y | ord_{i}u\leq r_{i}$ for
$i = 1,\dots, p+1$ and $u$ is not a multiple of any $v_{A_{j}}$, $j = 1,\dots, d\}$ and
$U_{r_{1}\dots r_{p+1}}^{(2)} = \{u\in \Lambda Y | ord_{i}u\leq r_{i}$ for
$i = 1,\dots, p+1$ and there exists at least one pair $i, j$ ($1\leq i\leq p,\,
1\leq j\leq d$) such that $u = \lambda v_{A_{j}}$, $\lambda\sim v_{A_{j}}$,  and
$ord_{i}(\lambda u_{A_{j}}^{(i)}) > r_{i}\}$.  Clearly, $U_{r_{1}\dots r_{p+1}} =
U_{r_{1}\dots r_{p+1}}^{(1)}\bigcup U_{r_{1}\dots r_{p+1}}^{(2)}$ and
$U_{r_{1}\dots r_{p+1}}^{(1)}\bigcap U_{r_{1}\dots r_{p+1}}^{(2)} = \emptyset$.

By Theorem 4.6, there exists a numerical polynomial
$\phi(t_{1},\dots, t_{p+1})$ in $p+1$ variables $t_{1},\dots, t_{p+1}$ such that
$\phi(r_{1},\dots, r_{p+1}) = Card\,U_{r_{1}\dots r_{p+1}}^{(1)}$ for all
sufficiently large $(r_{1}\dots r_{p+1})\in {\bf N}^{p+1}$,
$deg_{t_{i}}\phi \leq m_{i}$ for $i=1,\dots, p$, and $deg_{t_{p+1}}\phi\leq n$. Thus, in order to complete
the proof of the theorem, we need to show that there exists a numerical
polynomial $\psi(t_{1},\dots, t_{p+1})$  in $p+1$ variables
$t_{1},\dots, t_{p+1}$ such that $\psi(r_{1},\dots, r_{p+1}) =
Card\,U_{r_{1}\dots r_{p+1}}^{(2)}$
for all sufficiently large $(r_{1}\dots r_{p+1})\in {\bf N}^{p+1}$,
$deg_{t_{i}}\psi \leq m_{i}$ for $i=1,\dots, p$, and $deg_{t_{p+1}}\psi\leq n$.

Let $ord_{i}v_{A_{j}} = a_{ij}$,
$ord_{i}u_{A_{j}}^{(i)} = b_{ij}$, and $ord_{\sigma}v_{A_{j}} = c_{j}$ for $i=1,\dots, p; j= 1,\dots, d$
(clearly, $a_{ij}\leq b_{ij}$ for $i=1,\dots, p;\, j=1,\dots, d$). Furthermore, for any
$q=1,\dots, p$ and for any integers $k_{1}, \dots, k_{q}$ such that
$1\leq k_{1} < \dots < k_{q}\leq p$, let
$V_{j; k_{1},\dots, k_{q}}(r_{1},\dots, r_{p+1}) =
\{\lambda v_{A_{j}}|\lambda\sim v_{A_{j}}$, $ord_{i}\lambda \leq r_{i}-a_{ij}$ for $i=1,\dots, p$,
$ord_{\sigma}\lambda\leq r_{p+1}-c_{j}$,
and  $ord_{k}\lambda > r_{k}-b_{kj}$ if and only if $k$ is equal to one
of the numbers $k_{1}, \dots, k_{q}\}$.
Using Theorem 4.3(iii) we obtain that $Card\,V_{j; k_{1},\dots, k_{q}}(r_{1},\dots, r_{p+1}) =
\phi_{j; k_{1},\dots, k_{q}}(r_{1},\dots, r_{p+1})$, where
$\phi_{j; k_{1},\dots, k_{q}}(t_{1},\dots,$

\noindent$t_{p+1})$ is a numerical polynomial
in $p$ variables $t_{1},\dots, t_{p+1}$ defined by the formula
$$\phi_{j; k_{1},\dots, k_{q}}(t_{1},\dots, t_{p}) =
{{t_{1}+m_{1}-a_{1j}}\choose{m_{1}}}\dots
{{t_{k_{1} -1}+m_{k_{1} -1}-a_{k_{1} -1,j}}\choose{m_{k_{1} -1}}}$$

\noindent$$\left[{{t_{k_{1}}+m_{k_{1}}-a_{k_{1},j}}\choose{m_{k_{1}}}} -
{{t_{k_{1}}+m_{k_{1}}-b_{k_{1},j}}\choose{m_{k_{1}}}}\right]
{{t_{k_{1}+1}+m_{k_{1}+1}-a_{k_{1}+1,j}}\choose{m_{k_{1}+1}}}$$

$$\dots {{t_{k_{q} -1}+m_{k_{q} -1}-a_{k_{q} -1,j}}\choose{m_{k_{q} -1}}}
\left[{{t_{k_{q}}+m_{k_{q}}-a_{k_{q},j}}\choose{m_{k_{q}}}} -\\
{{t_{k_{q}}+m_{k_{q}}-b_{k_{q},j}}\choose{m_{k_{q}}}}\right]$$
\begin{equation}\dots {{t_{p}+m_{p}-a_{pj}}\choose{m_{p}}}{{t_{p+1}+n-c_{j}}\choose{n}}
\end{equation}
Clearly, $deg_{t_{i}}\phi_{j; k_{1},\dots, k_{q}} \leq m_{i}$
for $i=1, \dots, p$ and $deg_{t_{p+1}}\phi_{j; k_{1},\dots, k_{q}} \leq n$.

\medskip

Now, for any $j=1,\dots, d$, let
$V_{j}(r_{1}, \dots, r_{p+1}) = \{\lambda v_{A_{j}}|\lambda\sim v_{A_{j}}$,
$ord_{i}\lambda \leq r_{i}-a_{ij}$ for $i=1,\dots, p$,
$ord_{\sigma}\lambda\leq r_{p+1}-c_{j}$, and there exists
$k\in {\bf N}$, $1\leq k\leq p$, such that $ord_{k}\lambda > r_{k}-b_{kj}\}$.
Then the combinatorial principle of inclusion and exclusion implies
that $Card\,V_{j}(r_{1}, \dots, r_{p+1}) = \phi_{j}(r_{1}, \dots, r_{p+1})$,
where  $\phi_{j}(t_{1}, \dots, t_{p+1})$ is a numerical polynomial
in $p+1$ variables $t_{1},\dots, t_{p+1}$ defined by the formula

$$\phi_{j}(t_{1}, \dots, t_{p+1}) =
\D\sum_{k_{1}=1}^{p}\phi_{j;k_{1}}(t_{1}, \dots, t_{p+1}) -
\D\sum_{1\leq k_{1}<k_{2}\leq p}\phi_{j;k_{1},k_{2}}(t_{1}, \dots, t_{p+1})
+\dots +$$
\begin{equation}(-1)^{\nu-1}\D\sum_{1\leq k_{1}<\dots k_{\nu}\leq p}
\phi_{j;k_{1},\dots, k_{\nu}}(t_{1}, \dots, t_{p+1})+\dots +
(-1)^{p}\phi_{j;1,2,\dots,p}(t_{1}, \dots, t_{p+1}).
\end{equation}
It is easy to see that $deg_{t_{i}}\phi_{j}\leq m_{i}$
for $i=1, \dots, p$ and $deg_{t_{p+1}}\phi_{j} \leq n$.

Applying the principle of inclusion and exclusion we obtain that

$$Card\,U_{r_{1}\dots r_{p+1}}^{(2)} =
Card\,\bigcup_{j=1}^{d}V_{j}(r_{1}, \dots, r_{p+1}) =
\D\sum_{j=1}^{d}Card\,V_{j}(r_{1}, \dots, r_{p+1})$$
$$- \D\sum_{1\leq j_{1}<j_{2}\leq d}Card\,(V_{j_{1}}(r_{1}, \dots, r_{p+1})
\bigcap V_{j_{2}}(r_{1}, \dots, r_{p+1})) + \dots +$$
\begin{equation}(-1)^{d-1}Card\,\bigcap_{\nu=1}^{d}V_{\nu}(r_{1}, \dots, r_{p+1}),\hspace{1.1in}
\end{equation}
so it is sufficient to prove that for any $s = 1, \dots, d$ and for any
indices $j_{1}, \dots, j_{s},$\\$1\leq j_{1}< \dots < j_{s}\leq d$,
$Card\,(V_{j_{1}}(r_{1}, \dots, r_{p+1})\bigcap \dots
\bigcap V_{j_{s}}(r_{1}, \dots, r_{p+1})) =$

\noindent$\phi_{j_{1},\dots,j_{s}}(r_{1}, \dots, r_{p+1})$, where
$\phi_{j_{1},\dots,j_{s}}(t_{1}, \dots, t_{p+1})$ is a numerical polynomial
in $p+1$ variables $t_{1}, \dots, t_{p+1}$ such that
$deg_{t_{i}}\phi_{j_{1},\dots,j_{s}}\leq m_{i}$ for $i=1, \dots, p$ and
$deg_{t_{p+1}}\phi_{j_{1},\dots,j_{s}}\leq n$.

\medskip

It is clear that the intersection $V_{j_{1}}(r_{1}, \dots, r_{p+1})\bigcap \dots
\bigcap V_{j_{s}}(r_{1}, \dots, r_{p+1})$ is not empty (therefore,
$\phi_{j_{1},\dots,j_{s}} \neq 0$) if and only if
the leaders $v_{A_{j_{1}}}, \dots, v_{A_{j_{s}}}$ contain
the same $\Delta$-$\sigma$-indeterminate $y_{i}$ ($1\leq i\leq s$)
and they are all similar to each other. Let us consider
such an intersection and let $v(j_{1}, \dots, j_{s}) =
lcm(v_{A_{j_{1}}}, \dots, v_{A_{j_{s}}})$ (this least common multiple
of similar to each other terms is defined in a natural way, as the smallest with respect
to $<_{\sigma}$ common multiple of $v_{A_{j_{1}}}, \dots, v_{A_{j_{s}}}$).
Let elements $\gamma_{1}, \dots, \gamma_{s}\in \Lambda$ be defined
by the conditions $v(j_{1}, \dots, j_{s}) = \gamma_{\nu}v_{A_{j_{\nu}}}$ and
$\gamma_{\nu}\sim v_{A_{j_{\nu}}}$ ($\nu = 1, \dots, s$).

\medskip

Then  $V_{j_{1}}(r_{1}, \dots, r_{p+1})\bigcap \dots
\bigcap V_{j_{s}}(r_{1}, \dots, r_{p+1})$ is the set of all terms
$u = \lambda v(j_{1}, \dots, j_{s})$ such that $u\sim v(j_{1}, \dots, j_{s})$, $ord_{i}u\leq r_{i}$
(that is, $ord_{i}\lambda \leq r_{i} - ord_{i}v(j_{1}, \dots, j_{s})$)
for $i=1,\dots, p$, $ord_{\sigma}u\leq r_{p+1}$, and for any $l = 1, \dots, s$, there exists at least
one index $k\in \{1,\dots, p\}$ such that
$ord_{k}(\lambda\gamma_{l}u_{A_{j_{l}}}^{(k)}) > r_{k}$
(i. e., $ord_{k}\lambda > r_{k} - ord_{k}v(j_{1}, \dots, j_{s})
- ord_{k}u_{A_{j_{l}}}^{(k)} + ord_{k}v_{A_{j_{l}}}$).

Setting $c_{j_{1}, \dots, j_{s}}^{(i)} = ord_{i}v(j_{1}, \dots, j_{s})$
($1\leq i\leq p$), $c_{j_{1}, \dots, j_{s}}^{(p+1)} = ord_{\sigma}v(j_{1}, \dots, j_{s})$
and applying the principle of inclusion and exclusion
once again, we obtain that

\noindent$Card\,\bigcap_{\nu =1}^{s}V_{j_{\nu}}(r_{1}, \dots, r_{p+1})$
is an alternating sum of terms of the form

\noindent$Card\,W(j_{1}, \dots, j_{s}; k_{11}, k_{12},\dots, k_{1q_{1}}, k_{21},
\dots, k_{sq_{s}}; r_{1},\dots, r_{p+1})$ where\\
$W(j_{1}, \dots, j_{s}; k_{11}, k_{12},\dots, k_{1q_{1}}, k_{21},
\dots, k_{sq_{s}}; r_{1},\dots, r_{p+1})$ denotes the set

\smallskip

\noindent $\{\lambda \in \Lambda |\lambda\sim v(j_{1}, \dots, j_{s}),
ord_{i}\lambda \leq r_{i} - c_{j_{1}, \dots, j_{s}}^{(i)}$ for $i=1,\dots, p$,
$ord_{\sigma}\lambda \leq r_{p+1}-c_{j_{1}, \dots, j_{s}}^{(p+1)}$ and for any
$l=1, \dots, s, 1\leq k\leq p, \,ord_{k}\lambda > r_{k} -
c_{j_{1}, \dots, j_{s}}^{(k)} + a_{kj_{l}} - b_{kj_{l}}$
if and only if $k$ is equal to one of the numbers
$k_{l1}, \dots k_{lq_{l}}\}$

\smallskip

\noindent($q_{1}, \dots, q_{s}$ are some positive
integers in the set $\{1,\dots, p\}$
and $\{k_{i\nu} | 1\leq i\leq s, 1\leq \nu \leq q_{s}\}$ is a family of
integers such that $1\leq k_{i1} < k_{i2} < \dots < k_{iq_{i}}\leq p$
for $i=1,\dots, s$).

\noindent Thus, it remains to show that
$Card\,W(j_{1}, \dots, j_{s}; k_{11},\dots, k_{sq_{s}}; r_{1},\dots, r_{p+1})
=$

\noindent$\psi_{ k_{11}, \dots, k_{sq_{s}}}^{j_{1}, \dots, j_{s}}(r_{1}, \dots, r_{p+1})$
where  $\psi_{ k_{11}, \dots, k_{sq_{s}}}^{j_{1}, \dots, j_{s}}
(t_{1}, \dots, t_{p+1})$
is a numerical polynomial in $p+1$ variables $t_{1}, \dots, t_{p+1}$ whose degrees
with respect to $t_{i}$ ($1\leq i\leq p$) and $t_{p+1}$ do not
exceed $m_{i}$ and $n$, respectively.  But this is almost evident: as in the 
evaluation of  $Card\,V_{j; k_{1},\dots, k_{q}}(r_{1},\dots, r_{p+1})$
(when we use Theorem 4.3 (iii) to get formula (5.3)), we see that
$Card\,W(j_{1}, \dots, j_{s}; k_{11},\dots, k_{sq_{s}}; r_{1},\dots, r_{p})$
is a product of terms of the form $\D{{r_{p+1}+n-c_{j_{1},
\dots, j_{s}}^{(p+1)}}\choose{n}}$ or $\D{{r_{k}+m_{k}-c_{j_{1},
\dots, j_{s}}^{(k)}-S_{k}}\choose{m_{k}}}$
(the last term corresponds to an integer $k$ such that $1\leq k\leq p$ and
$k\neq k_{i\nu}$ for any $i=1,\dots, s, \nu = 1,\dots, q_{s}$; the number
$S_{k}$ is  $\max \{b_{kj_{l}} - a_{kj_{l}}| 1\leq l\leq s\}$)
or
$\left[\D{{r_{k}+m_{k}-c_{j_{1},\dots, j_{s}}^{(k)}}\choose{m_{k}}} -
\D{{r_{k}+m_{k}-c_{j_{1},\dots, j_{s}}^{(k)}-T_{k}}\choose{m_{k}}}\right]$
(such a term appears in the product if $k$ is equal to some $k_{i\nu}$
($1\leq i\leq s, 1\leq \nu\leq q_{s}$). In this case, if
$k_{i_{1}\nu_{1}}, \dots, k_{i_{e}\nu_{e}}$ are all elements of the
set $\{k_{i\nu} | 1\leq i\leq s, 1\leq \nu\leq q_{s}\}$ that are equal
to $k$ ($1\leq e\leq s, 1\leq i_{1} < \dots < i_{e}\leq s$), then
$T_{k}$ is defined as $\min \{b_{kj_{i_{\lambda}}} - a_{kj_{i_{\lambda}}} |
1\leq \lambda \leq l\}$).

The corresponding numerical polynomial
$\psi_{ k_{11}, \dots, k_{sq_{s}}}^{j_{1}, \dots, j_{s}}
(t_{1}, \dots, t_{p+1})$ is a product of $p$ ''elementary'' numerical
polynomials $f_{1}, \dots, f_{p}$ where $f_{k}$ ($1\leq k\leq p$) is 
a polynomial of the form $\D{{t_{p+1}+n-c_{j_{1},
\dots, j_{s}}^{(p+1)}}\choose{n}}$ or
$\D{{t_{k}+m_{k}-c_{j_{1}, \dots, j_{s}}^{(k)}-S_{k}}\choose{m_{k}}}$  or

\medskip

\noindent$\left[\D{{t_{k}+m_{k}-c_{j_{1},\dots, j_{s}}^{(k)}}\choose{m_{k}}} -
\D{{t_{k}+m_{k}-c_{j_{1},\dots, j_{s}}^{(k)}-T_{k}}\choose{m_{k}}}\right]$ ($1\leq k\leq p$).
Since the degree of such a product with respect to any
variable $t_{i}$ ($1\leq i\leq p$) and $t_{p+1}$ does not exceed $m_{i}$ and $n$,
respectively, this completes the proof of the first two parts of Theorem 3.1.

\medskip

In order to prove the last part of the theorem, suppose that $\zeta = \{\zeta_{1},\dots, \zeta_{q}\}$ is
another system of $\Delta$-$\sigma$-generators of $L/K$, that is,
$L = K\langle \eta_{1},\dots,\eta_{s}\rangle =
K\langle \zeta_{1},\dots,\zeta_{q}\rangle$). Let
$$\Phi_{\zeta}(t_{1},\dots, t_{p+1}) =
\D\sum_{i_{1}=0}^{m_{1}}\dots \D\sum_{i_{p}=0}^{m_{p}}b_{i_{1}\dots i_{p+1}}
{t_{1}+i_{1}\choose i_{1}}\dots {t_{p+1}+i_{p+1}\choose i_{p+1}}$$
be the dimension polynomial of our $\Delta$-$\sigma$-field  extension associated
with the system of generators $\zeta$. Then there exist
$h_{1}, \dots, h_{p+1}\in {\bf N}$ such that
$\eta_{i} \in K(\bigcup_{j=1}^{q}\Lambda(h_{1},\dots, h_{p+1})\zeta_{j})$
and $\zeta_{k} \in K(\bigcup_{j=1}^{s}\Lambda(h_{1},\dots, h_{p+1})\eta_{j})$
for any $i=1,\dots, s$ and $k=1,\dots, q$, whence
$\Phi_{\eta}(r_{1},\dots, r_{p+1}) \leq
 \Phi_{\zeta}(r_{1}+h_{1},\dots, r_{p+1}+h_{p+1})$ and
$\Phi_{\zeta}(r_{1},\dots, r_{p+1}) \leq
 \Phi_{\eta}(r_{1}+h_{1},\dots, \\r_{p+1}+h_{p+1})$ for all sufficiently large
$(r_{1},\dots, r_{p+1})\in {\bf N}^{p+1}$. Now the statement of the third
part of Theorem 3.1 follows from the fact that for any element
$(k_{1},\dots, k_{p+1})\in E_{\eta}'$, the term
$\D{t_{1}+k_{1}\choose k_{1}}\dots {t_{p+1}+k_{p+1}\choose k_{p+1}}$ appears in
$\Phi_{\eta}(t_{1},\dots, t_{p+1})$ and $\Phi_{\zeta}(t_{1},\dots, t_{p+1})$
with the same coefficient $a_{k_{1}\dots k_{p+1}}$.

\bigskip

Let us consider some examples of computation of $\Delta$-$\sigma$-dimension polynomials.

\begin{example}
{\em Consider a $\Delta$-$\sigma$-field extension $L = K\langle \eta\rangle$ where $\Delta = \{\delta\}$,
$\sigma = \{\alpha\}$, and the $\Delta$-$\sigma$-generator of $L/K$ satisfies the defining equation
\begin{equation}
\delta\eta - \alpha\eta  - a = 0
\end{equation}
where $a\in K$.

Since the equation is linear, the corresponding defining  $\Delta$-$\sigma$-ideal of the ring of
$\Delta$-$\sigma$-polynomials $K\{y\}$ is $P = [\alpha y - \delta y - a]$ (the fact that any linear
ideal is prime is well-known, see, for example, \cite[Chapter IV, Section 5]{K2}). By Proposition 5.10,
the $\Delta$-$\sigma$-polynomials $A = \alpha y - \delta y - a$ and $-\alpha^{-1}A = \alpha^{-1}\delta y - y - \alpha^{-1}(a)$
form a characteristic set of the ideal $P$. Using the notation of the proof of Theorem 3.1 and applying the procedure
described in this proof together with Theorem 4.3(iii), and formula (4.4),
we obtain that $Card\,U_{r_{1}r_{2}}^{(1)}  = r_{1} + r_{2} +1$ and $Card\,U_{r_{1}r_{2}}^{(2)} = r_{2}$ for all sufficiently large
$(r_{1}, r_{2})\in {\bf N}^{2}$. Therefore the corresponding $\Delta$-$\sigma$-dimension polynomial is as follows:
$\Phi_{\eta}(t_{1}, t_{2}) = t_{1} + 2t_{2} + 1$. Note that this polynomial expresses the strength of the difference-differential equation
${\D\frac{dy(x)}{dx}} - y(x+h) - a(x) = 0$ where $y(x)$ is an unknown function and $a(x)$ belongs to a functional $\Delta$-$\sigma$-field
$K$ with the derivation ${\D\frac{d}{dx}}$ and automorphism $f(x)\mapsto f(x+h)$ where $h$ is a constant of the field (say, if $K$ is a field of
functions of a real variable, then $h$ is a real number).}
\end{example}

\begin{example}
{\em Let us find the $\Delta$-$\sigma$-dimension polynomial that expresses the strength of the difference-differential equation
\begin{equation}
\frac{\partial^{2} y(x_{1}, x_{2})}{\partial x_{1}^{2}} + \frac{\partial^{2} y(x_{1}, x_{2})}{\partial x_{2}^{2}} + y(x_{1} + h) + a(x)  = 0
\end{equation}
over some $\Delta$-$\sigma$-field of functions of two real variables $K$, where the basic set of derivations
$\Delta = \{\delta_{1} = \frac{\partial}{\partial x_{1}}, \delta_{2} = \frac{\partial}{\partial x_{2}}\}$ has the partition
$\Delta = \{\delta_{1}\}\bigcup \{\delta_{1}\}$ and $\sigma$ consists of one automorphisms $\alpha: f(x_{1}, x_{2})\mapsto f(x_{1}+h, x_{2})\}$
($\alpha$ is the shift of the first argument of a function by a real number $h$).

In this case, the associated $\Delta$-$\sigma$-extension $K\langle\eta\rangle/K$ is $\Delta$-$\sigma$-isomorphic to the field of fractions of the
integral $\Delta$-$\sigma$-domain $K\{y\}/[\alpha y + \delta_{1}^{2}y + \delta_{1}^{2}y + a]$ (the element $a\in K$ corresponds to the function $a(x)$).
Applying Proposition 5.10 we obtain that the characteristic set of the defining ideal of the corresponding $\Delta$-$\sigma$-extension
$K\langle\eta\rangle/K$ consists of the $\Delta$-$\sigma$-polynomials $g_{1} = \alpha y + \delta_{1}^{2}y + \delta_{1}^{2}y + a$ and
$g_{2} = \alpha^{-1}g_{1} = \alpha^{-1}\delta_{1}^{2}y + \alpha^{-1}\delta_{2}^{2}y + y + \alpha^{-1}(a)$.
With the notation of the proof of Theorem 3.1, the application of the procedure described in this proof, Theorem 4.3(iii), and formula (4.4)
leads to the following expressions for the numbers of elements of the sets $U_{r_{1}r_{2}r_{3}}^{(1)}$ and $U_{r_{1}r_{2}r_{3}}^{(2)}$:
$Card\,U_{r_{1}r_{2}r_{3}}^{(1)} = r_{1}r_{2}  + 2r_{2}r_{3} + r_{1} + r_{2} + 2r_{3} + 1$
and $Card\,U_{r_{1}r_{2}r_{3}}^{(2)} = 4r_{1}r_{3}  + 2r_{2}r_{3} - 2r_{3}$ for all sufficiently large $(r_{1}, r_{2}, r_{3})\in {\bf N}^{3}$.
Thus, the strength of equation (5.6) corresponding to the given partition of the basic set of derivations is expressed by the $\Delta$-$\sigma$-polynomial
$\Phi_{\eta}(t_{1}, t_{2}, t_{3}) = t_{1}t_{2} + 4t_{1}t_{3} + 4t_{2}t_{3} + t_{1} + t_{2} + 1$.}
\end{example}

\section{Generalized Gr\"obner bases in free difference-differential modules and difference-differential dimension
polynomials}

\setcounter{equation}{0}

Let $K$ be a difference-differential field of zero characteristic with basic
sets $\Delta = \{\delta_{1},\dots, \delta_{m}\}$ and $\sigma = \{\alpha_{1},\dots, \alpha_{n}\}$
of derivations and automorphisms, respectively. Suppose that partition (3.1) of the set of derivations is
fixed:

\noindent$\Delta = \Delta_{1}\bigcup \dots \bigcup \Delta_{p}$ \,\,\, ($p\geq 1$)
where $\Delta_{1} = \{\delta_{1},\dots, \delta_{m_{1}}\}, \,\Delta_{2} = \{\delta_{m_{1}+1},\dots,$

\noindent$\delta_{m_{1}+m_{2}}\},\dots, \Delta_{p} = \{\delta_{m_{1}+\dots + m_{p-1}+1},\dots, \delta_{m}\}$\,
($m_{1}+\dots + m_{p} = m$).

As before, let $\Lambda$ denote the free commutative semigroup of all power products of the form
$\lambda = \delta_{1}^{k_{1}}\dots \delta_{m}^{k_{m}}\alpha_{1}^{l_{1}}\dots \alpha_{n}^{l_{n}}$ where
$k_{i}\in {\bf N},\, l_{j}\in {\bf Z}$ ($1\leq i\leq m,\, 1\leq j\leq n$) and let $ord_{i}\lambda$ ($1\leq i\leq p$)
and $ord_{\sigma}\lambda$ denote the orders of such an element $\lambda$ relative to $\Delta_{i}$ and $\sigma$, respectively,
introduced at the beginning of section 3.

In what follows an expression of the form
$\sum_{\lambda \in \Lambda}a_{\lambda}\lambda$, where $a_{\lambda}\in K$
for all $\lambda \in \Lambda$ and only finitely many coefficients
$a_{\lambda}$ are different from zero, is called a difference-differential
(or $\Delta$-$\sigma$-) operator over $K$. Two $\Delta$-$\sigma$-operators
$\sum_{\lambda \in \Lambda}a_{\lambda}\lambda$ and $\sum_{\lambda \in \Lambda}b_{\lambda}\lambda$
are considered to be equal if and only if $a_{\lambda} = b_{\lambda}$ for all
$\lambda \in \Lambda$.

The set of all $\Delta$-$\sigma$-operators over $K$ can be equipped with a ring structure if one sets
$\sum_{\lambda \in \Lambda}a_{\lambda}\lambda +
\sum_{\lambda \in \Lambda}b_{\lambda}\lambda =
\sum_{\lambda \in \Lambda}(a_{\lambda} + b_{\lambda})\lambda$,
$a(\sum_{\lambda \in \Lambda}a_{\lambda}\lambda) =
\sum_{\lambda \in \Lambda}(aa_{\lambda})\lambda$,
$(\sum_{\lambda \in \Lambda}a_{\lambda}\lambda)\mu =
\sum_{\lambda \in \Lambda}a_{\lambda}(\lambda \mu)$,
$\mu a = \mu(a)\mu$ for any $\Delta$-$\sigma$-operators
$\sum_{\lambda \in \Lambda}a_{\lambda}\lambda$,
$\sum_{\lambda \in \Lambda}b_{\lambda}\lambda$ and for any elements
$a\in K, \mu \in \Lambda$, and extend these rules by the distributivity.
This ring is called the ring of difference-differential
(or $\Delta$-$\sigma$-) operators over $K$, it will be denoted by $D$
(or $D_{K}$, if the ring $K$ should be specified).

If $u = \sum_{\lambda \in \Lambda}a_{\lambda}\lambda$ is a
$\Delta$-$\sigma$-operator over $K$, then the orders of $u$ relative to
the sets $\Delta_{i}$ ($1\leq i\leq p$) and $\sigma$ are defined as numbers
$ord_{i}u = max\{ord_{i}\lambda | a_{\lambda}\neq 0\}$ and
$ord_{\sigma}u = max\{ord_{\sigma}\lambda | a_{\lambda}\neq 0\}$,
respectively. The number $ord\,u = max\{ord\,\lambda | a_{\lambda}\neq 0\}$ is said
to be the order of the $\Delta$-$\sigma$-operator $u$.

For any $r_{1},\dots, r_{p+1}\in {\bf N}$, let $D_{r_{1}\dots r_{p+1}}$ denote
the vector $K$-subspace of $D$ generated by $\Lambda(r_{1},\dots,
r_{p+1})$. Setting $D_{r_{1}\dots r_{p+1}} = 0$ for $(r_{1},\dots,
r_{p+1})\in {\bf Z}^{p+1}\setminus {\bf N}^{p+1}$, we obtain a family
$\{D_{r_{1}\dots r_{p+1}} | (r_{1},\dots, r_{p+1})\in {\bf
Z}^{p+1}\}$ called a {\em standard $(p+1)$-dimensional filtration\/} of
$D$.

Obviously, $\bigcup{\{D_{r_{1}\dots r_{p+1}}|r_{1}\dots r_{p+1}\in {\bf Z}\}}=D$,
$D_{r_{1}\dots r_{p+1}} \subseteq D_{r_{1}\dots r_{i-1},r_{i}+1,r_{i+1},\dots r_{p+1}}$
for all $i$, $1\leq i\leq, p+1$,  and $D_{r_{1}\dots r_{p+1}}D_{s_{1}\dots s_{p+1}}=
D_{r_{1}+s_{1},\dots, r_{p+1}+s_{p+1}}$ for any $r_{i}, s_{i}\in {\bf N},\, 1\leq i\leq p+1$.

\medskip

A left module over the ring $D$ is called a {\em difference-differential $K$-module}
or a {\em $\Delta$-$\sigma$-$K$-module}.  In other words, a vector $K$-space
$M$ is called a $\Delta$-$\sigma$-$K$-module, if the elements of the set
$\Delta \bigcup \sigma ^{\ast}$ act on $M$ in such a way that
$\beta(x+y) = \beta(x) + \beta(y)$, $\beta(\gamma(x)) = \gamma(\beta(x))$,
$\delta(ax) = a\delta(x) + \delta(a)x$, $\tau(ax) = \tau(a)\tau(x)$, and
$\tau(\tau ^{-1}(x)) = x$ for any
$\beta, \gamma \in \Delta \bigcup \sigma ^{\ast}$,
$\delta \in \Delta$, $\tau \in \ \sigma ^{\ast}$, $a\in K$, and $x\in M$.

\begin{definition}
If  $M$ is a $\Delta$-$\sigma$-$K$-module, then a family $\{M_{r_{1}\dots r_{p+1}}|(r_{1}\dots r_{p+1})\in {\bf Z}\}$
of vector $K$-subspaces of the module $M$ is
called a  $(p+1)$-dimensional filtration of $M$ if the following three conditions hold:

(i) For any fixed integers $r_{1},\dots, r_{i-1}$,
$r_{i+1},\dots, r_{p+1}$\, ($1\leq i\leq p+1$), $M_{r_{1}\dots
r_{i}\dots r_{p+1}}\subseteq M_{r_{1}\dots r_{i-1}, r_{i}+1,
r_{i+1}\dots r_{p+1}}$ and $M_{r_{1}\dots r_{p+1}} = 0$ for all
sufficiently small $r_{i}\in {\bf Z}$;

\medskip

(ii) $\bigcup \{M_{r_{1}\dots r_{p+1}} | (r_{1},\dots, r_{p+1})
\in {\bf Z}^{p+1}\} = M$;

\medskip

 (iii) $D_{r_{1}\dots r_{p+1}}M_{s_{1}\dots s_{p+1}} \subseteq
M_{r_{1}+s_{1},\dots, r_{p+1}+ s_{p+1}}$ for any $(r_{1},\dots,
r_{p+1}), $ 

\noindent$(s_{1},\dots, s_{p+1})\in {\bf Z}^{p+1}$.

If every vector $K$-space $M_{r_{1}\dots r_{p+1}}$ is finitely
generated and there exists an element $(h_{1},\dots, h_{p+1})\in{\bf
Z}^{p+1}$ such that $D_{r_{1}\dots r_{p+1}}M_{h_{1}\dots h_{p+1}} =
M_{r_{1}+h_{1},\dots, r_{p+1}+ h_{p+1}}$ for any \\$(r_{1},\dots,
r_{p+1})\in{\bf N}^{p+1}$, then the $(p+1)$-dimensional filtration is
called excellent.
\end{definition}

It is easy to see that if $u_{1},\dots, u_{q}$ is a finite system
of generators of a left $D$-module $M$, then the filtration
$\{\sum_{i=1}^{q}D_{r_{1}\dots r_{p+1}}u_{i} |$ $(r_{1},\dots,
r_{p+1})\in {\bf Z}^{p+1}\}$ is excellent.

\medskip

Let $M$ and $N$ be two $\Delta$-$\sigma$-$K$-modules.
A homomorphism of of vector $K$-spaces $f: M \longrightarrow N$ is called a
{\em $\Delta$-$\sigma$-homomorphism} (or {\em difference-differential homomorphism}),
if $f(\tau x) = \tau f(x)$ for any
$x\in M, \beta \in \Delta \bigcup \sigma ^{\ast}$. Surjective (respectively,
injective or bijective) $\Delta$-$\sigma$-homomorphism is called a
$\Delta$-$\sigma$-epimorphism ($\Delta$-$\sigma$-monomorphism or
$\Delta$-$\sigma$-isomorphism, respectively).

\medskip

Let $F$ be a finitely generated free left $D$-module with free generators $f_{1},\dots, f_{q}$.
(Using the ''difference-differential'' terminology, we also say that $F$ is a free difference-differential $K$-module
(or a free $\Delta$-$\sigma$-$K$-module) with the set of free $\Delta$-$\sigma$-generators $\{f_{1},\dots, f_{q}\}$.)

Then $F$ can be considered as a vector $K$-space with the basis
$\{\lambda f_{i}\,|\,\lambda\in \Lambda$, $1\leq i\leq q\}$).  This set will be denoted
by $\Lambda f$ and its elements will be called {\em terms\/}. (The fact that the same name is used in section 5 for certain elements of
the ring of $\Delta$-$\sigma$-polynomials will not cause any confusions.) For any
term $\lambda f_{i}$, the element $\lambda$ will be called the {\em head} of
the term.

Since $\Lambda f$ is a basis of $F$ over $K$, every element $f\in F$ has a unique representation
as a linear combination of terms:
\begin{equation}
f = a_{1}\lambda_{1}f_{i_{1}} + \dots + a_{d}\lambda_{1}f_{i_{d}}
\end{equation}
for some nonzero elements $a_{j}\in K$ and some $\lambda_{j}\in \Lambda$
($1\leq j\leq d$). We say that the element $f$ contains a term
$\lambda_{k}f_{k}$, if this term appears in the representation (6.1)
with nonzero coefficient.

We define the orders of a term $\lambda f_{j}\in \Lambda f$ relative to
the sets $\Delta_{i}$  and $\sigma$ as the corresponding orders of
$\lambda$: $ord_{i}(\lambda f_{j}) = ord_{i}\lambda$ ($1\leq i\leq p$) and
$ord_{\sigma}(\lambda f_{j}) = ord_{\sigma}\lambda$. The number
$ord\,(\lambda f_{j}) = ord\,\lambda$  is said to be the order of the term
$\lambda f_{j}$.

We say that two terms $u=\lambda f_{i}$ and $v=\lambda' f_{j}$
are {\em similar} and write $u\sim v$ if $\lambda\sim\lambda'$.
If $u = \lambda e_{i}$ is a term and $\lambda'\in \Lambda$, we say that $u$ is similar to $\lambda'$ and write
$u\sim \lambda'$ if $\lambda\sim \lambda'$. Furthermore, if $u, v\in \Lambda f$, we say that $u$ {\em divides} $v$ or
{\em $v$ is a multiple of $u$}, if $u=\lambda e_{i}$, $v=\lambda' e_{i}$ for some $e_{i}$ and $\lambda|\lambda'$
(in particular, it means that $\lambda\sim \lambda'$).

In what follows we will consider $p+1$ orders $<_{1},\dots, <_{p}, <_{\sigma}$ on the set  $\Lambda f$ defined in the same way as the
corresponding orders of the terms in the ring of $\Delta$-$\sigma$-polynomials: $\lambda f_{j} <_{i}$(or $<_{\sigma}$) $\lambda' f_{k}$
($\lambda, \lambda'\in \Lambda,\, 1\leq j, k\leq q$) if and only if
$\lambda <_{i}$ (respectively, $<_{\sigma}$) $\lambda'$ in $\Lambda$, or $\lambda = \lambda'$ and $j < k$.

If an element $f\in F$ is written in the form (6.1) then the greatest element of the set
$\{\lambda_{1}f_{i_{1}}, \dots, \lambda_{d}f_{i_{d}}\}$ relative
to the orders $<_{i}$ ($1\leq i\leq p$) and $<_{\sigma}$, are called the $\Delta_{i}$-leader
and $\sigma$-leader of $f$; they are denoted by $u_{f}^{(i)}$ and $v_{f}$, respectively.

\begin{definition}
Let $f$ and $g$ be two elements of the free $\Delta$-$\sigma$-$K$-module $F$.
The element $f$ is said to be reduced with respect to $g$ if $f$ does not contain any multiple $\lambda v_g$
($\lambda \in \Lambda$) of the $\sigma$-leader $v_g$ such that
$ord_{i}(\lambda u_{g})\leq ord_{i}u_f$ for $i=1,\dots p$. \, An element
$h\in E$ is said to be reduced with respect to a set $\Sigma \subseteq E$,
if $h$ is reduced with respect to every element of the set $\Sigma$.
\end{definition}
Let us consider $p$ new symbols $z_{1},\dots, z_{p}$, the
free commutative semigroup $\Gamma$ of all power products $\lambda z_{1}^{k_{1}}\dots
z_{p}^{k_{p}}$ ($\lambda\in \Lambda, \,k_{i}\in {\bf N}$ for $i=1,\dots, p$), and the set
$\Gamma f$ of all elements of the form  $\{\gamma f_{j} | \gamma \in
\Gamma, 1\leq j\leq q\} = \Gamma \times \{e_{1},\dots, e_{q}\}$. If $a = \lambda z_{1}^{k_{1}}\dots
z_{p}^{k_{p}}f_{\mu}, b = \lambda' z_{1}^{l_{1}}\dots z_{p}^{l_{p}}f_{\nu}\in \Gamma f$ ($1\leq \mu, \nu\leq q$),
we say that $a$ divides $b$ (or $b$ is a multiple of $a$) and write $a|b$ if $\mu=\nu$, $\lambda|\lambda'$ in $\Lambda$,
and $k_{i}\leq l_{i}$ for $i=1,\dots, p$.

For any $i=1,\dots, p$ and for any $f\in F$, let $d(f) = ord_{i}u_{f}^{(i)} - ord_{i}v_{f}$, and let a mapping
$\rho:F\rightarrow \Gamma f$ be defined by $\rho(f) = z_{1}^{d_{1}(f)}\dots z_{p}^{d_{p}(f)}v_{f}$.
\begin{definition}
With the above notation, let $N$ be a $D$-submodule of $F$. A
finite set $G = \{g_{1},\dots, g_{r}\}\subseteq N$ is called a $\Delta$-$\sigma$-Gr\"obner basis of $N$
(with respect to the given partition (3.1) of the basic set $\Delta$) if for any $f\in N$, there exists $g_{i}\in G$ such that
$\rho(g_{i}) | \rho(f)$.
\end{definition}
\begin{definition}
Given $f, g, h\in F$, with $g\neq 0$, we say that the element $f$ reduces to
$h$  modulo $g$ in one step and write $f\xrightarrow {\text{g}} h$
 if and only if $f$ contains some term $w$ with a
coefficient $a$ such that $v_{g}|w$, $w = \lambda v_{g}$ ($\lambda\in \Lambda,\,\lambda\sim v_{g}$),
$$h = f - a\left(\lambda_{\sigma}(lc_{\sigma}(g))\right)^{-1}\lambda g,$$
and \, $ord_{i}(\lambda u_{g}^{(i)})\leq ord_{i}u_{f}^{(i)}$ for $i=1,\dots, p$.
\end{definition}
\begin{definition}
Let $f, h\in F$ and let $G = \{g_{1},\dots, g_{r}\}$ be a finite
set of nonzero elements of $F$. We say that $f$ reduces to $h$ modulo $G$
and write $f\xrightarrow {\text{$G$}} h$ if and only if there exists a sequence of elements
$g^{(1)}, g^{(2)},\dots g^{(q)}\in G$ and a sequence of elements
$h_{1},\dots, h_{q-1}\in F$ such that
$f\xrightarrow {\text{$g^{(1)}$}} h_{1}\xrightarrow {\text{$g^{(2)}$}} \dots \xrightarrow {\text{$g^{(q-1)}$}}h_{q-1}\xrightarrow {\text{$g^{(q)}$}}h$.
\end{definition}

The introduced notion of reduction is similar to the concept of
$(<_{k}, <_{i_{1}}, \dots <_{i_{l}})$-reduction defined in \cite[Section 3]{Levin3}.
The following two propositions can be proved in the same way as Theorem 3.6 and Proposition 3.7 of \cite{Levin3}.
\begin{proposition}
With the above notation, let $G =\{g_{1},\dots,g_{r}\}\subseteq F$
be a $\Delta$-$\sigma$-Gr\"obner basis of a $D$-submodule $N$ of $F$ and let $f\in N$.  Then there exist elements $g\in N$ and
$Q_{1},\dots,Q_{r}\in D$ such that $f-g=\sum_{i=1}^{r}Q_{i}g_{i}$
and $g$ is reduced with respect to $G$.
\end{proposition}

Now we can propose an approach to the computation of the dimension polynomial of
a finitely generated $\Delta$-$\sigma$-field extension
$L = K\langle\eta_{1},\dots,\eta_{s}\rangle$, which is based on the
consideration of the corresponding module of K\"ahler differentials
$\Omega_{L|K}$. As it is shown in \cite{Johnson} and \cite[Proposition 3.4.46]{KLMP},
the $L$-vector space $\Omega_{L|K}$ can be equipped with the structure of a left module
over the ring of $\Delta$-$\sigma$-operators $\mathcal{D} = D_{L}$ over $L$ in such a way that
the corresponding action of the elements of $\Delta\bigcup\sigma$ on $\Omega_{L|K}$ satisfy the conditions
$\delta(d\zeta) = d\delta(\zeta)$ and $\alpha(d\zeta) = d\alpha(\zeta)$
for any $\zeta\in L$, $\delta\in \Delta$, $\alpha\in \sigma$.

Let us consider $\mathcal{D}$ as a filtered ring with a natural $(p+1)$-dimensional
filtration $\{\mathcal{D}_{r_{1}\dots r_{p+1}} | (r_{1},\dots, r_{p+1})\in {\bf Z}^{p+1}\}$ where
$\mathcal{D}_{r_{1}\dots r_{p+1}}$ is the vector $L$-subspace of
$\mathcal{D}$ generated by $\Lambda(r_{1},\dots, r_{p+1})$ if all
$r_{i}\geq 0$, and  $0$ otherwise. Then the finitely generated  $\mathcal{D}$-module
$M = \Omega_{L|K} = \sum_{i=1}^{s}\mathcal{D}d\eta_{i}$ has a natural
 $(p+1)$-dimensional filtration $\{M_{r_{1}\dots r_{p+1}}
|(r_{1},\dots,\\ r_{p+1})\in {\bf Z}^{p+1}\}$ where $M_{r_{1}\dots r_{p+1}}$ ($r_{i}\geq 0$)
is the vector $L$-subspace of $\Omega_{L|K}$ generated by the
set $\{d\eta| \eta \in K(\{\lambda\eta_{j}\,|\,\lambda\in
\Lambda(r_{1},\dots, r_{p+1}), 1\leq j\leq s\})$ and
$M_{r_{1}\dots r_{p+1}} = 0$ whenever at least one $r_{i}$
is negative. Clearly, this is an excellent filtration of $M$ and,
as it is shown in \cite{Johnson} and \cite[Theorem 6.4.11]{KLMP},
$$trdeg_{K}K(\D\bigcup_{j=1}^{s} \Lambda(r_{1},\dots, r_{p+1})\eta_{j}) =
dim_{L}(M_{r_{1}\dots r_{p+1}})$$ for all $(r_{1}\dots r_{p+1})\in {\bf N}^{p+1}$.

It follows that the computation of the difference-differential dimension polynomial of
the extension $L/K$ can be reduced to the computation of the dimension polynomial associated
with the filtration $\{M_{r_{1}\dots r_{p+1}}\}$; the existence of such a polynomial is established by
the following theorem whose proof is similar to the proof of Theorem 4.2 of \cite{Levin3}.

\bigskip

\begin{theorem} Let $\{M_{r_{1}\dots r_{p+1}} = \D\sum_{i=1}^{s}{\mathcal{D}}_{r_{1}\dots r_{p}}f_{i}\}$
be an excellent $(p+1)$-dimensional filtration of a left
$\mathcal{D}$-module $M=\D\sum_{i=1}^{s}{\mathcal{D}}f_{i}$ ($f_{1},\dots, f_{s}$ generate
$M$ as a ${\mathcal{D}}$-module). Then there exists a polynomial
$\phi(t_{1}, \dots, t_{p+1})\in{\bf Q}[t_{1}, \dots, t_{p+1}]$ such that

\bigskip

{\em (i)} $\phi(r_{1}, \dots, r_{p+1}) = dim_{L} M_{r_{1} \dots r_{p+1}}$ for
all sufficiently large $(r_{1},\dots, r_{p+1})\in {\bf Z}^{p}$.

\bigskip

{\em (ii)} \, $deg_{t_{i}}\phi \leq m_{i}$ ($1\leq i\leq p$) and
$deg_{t_{p+1}}\phi \leq n$, so
that $deg\,\phi \leq m+n$ and  $\phi(t_{1},\dots,
t_{p+1})$ can be represented as $$\phi = \sum_{i_{1}=0}^{m_{1}}\dots
\sum_{i_{p}=0}^{m_{p}}\sum_{i_{p+1}=0}^{n}a_{i_{1}\dots i_{p+1}} {t_{1}+i_{1}\choose
i_{1}}\dots \D{t_{p+1}+i_{p+1}\choose i_{p+1}}$$ where $a_{i_{1}\dots i_{p+1}}\in {\bf Z}$
and $2^{n}\,|\,a_{m_{1}\dots m_{p}n}.$

\medskip

{\em (iii)} \,Let $A = \{(i_{1},\dots, i_{p+1})\in {\bf N}^{p+1}\,|\,
0\leq i_{k}\leq m_{k}$ \,($1\leq k\leq p$),\, $0\leq i_{p+1}\leq n$
and $a_{i_{1}\dots i_{p+1}}\neq 0\}$. Then $d = deg\,\phi$, $a_{m_{1}\dots m_{p}n}$,
elements $(k_{1},\dots, k_{p+1})\in A'$, the corresponding
coefficients $a_{k_{1}\dots k_{p+1}}$ and the coefficients of the
terms of total degree $d$ do not depend on the choice of the excellent filtration.
Also, $\D\frac{a_{m_{1}\dots m_{p}n}}{2^{n}}$  is equal
to the maximal number of elements of $M$ linearly independent over
$\mathcal{D}$.

{\em (iv)} Let $E$ be a free left ${\mathcal{D}}$-module with free generators $e_{1}, \dots,
e_{s}$. Let $\pi: E\longrightarrow M$ be the natural
${\mathcal{D}}$-epimorphism ($\pi(e_{i}) = f_{i}$ for $i=1, \dots, s$), $N =
Ker\,\pi$, and $G = \{g_{1}, \dots, g_{d}\}$ a $\Delta$-$\sigma$-Gr\"obner basis of
$N$. Furthermore, for any
$(r_{1},\dots, r_{p+1})\in {\bf N}^{p+1}$, let
$U_{r_{1}\dots r_{p+1}}^{(1)} = \{u = \lambda e_{k}\,|\,ord_{i}\lambda\leq r_{i}$ for
$i = 1,\dots, p$, $ord_{\sigma}u\leq r_{p+1}$, and $u$ is not a multiple of
any $v_{g_{j}}$, $j = 1,\dots, d$ ($1\leq k\leq d)\}$ and let
$U^{(2)}_{r_{1}\dots r_{p+1}} = \{u = \lambda e_{k}\,|\,ord_{i}u\leq r_{i}$ for $i = 1,\dots, p$,
$ord_{\sigma}u\leq r_{p+1}$, and for every $\lambda\in\Lambda, g\in G$ such that $u = \lambda
u_{g}^{(1)}$, $\lambda\sim u_{g}^{(1)}$, there exists $i\in \{1,\dots, p\}$ such that
$ord_{i}\lambda u_{g}^{(i)}> r_{i}\}$.

Then for any $(r_{1}, \dots, r_{p+1})\in {\bf N}^{p+1}$, the
set $\pi(U^{(1)}_{r_{1}\dots r_{p+1}}\bigcup U^{(2)}_{r_{1}\dots r_{p+1}})$ is a basis of $M_{r_{1}\dots
r_{p+1}}$ over $K$. Therefore
$\phi(r_{1}, \dots, r_{p+1}) = Card\left(U^{(1)}_{r_{1}\dots r_{p+1}}\bigcup U^{(2)}_{r_{1}\dots r_{p+1}}\right)$ for all
sufficiently large $(r_{1},\dots, r_{p+1})\in {\bf Z}^{p}$.
\end{theorem}

The last theorem shows that the polynomial $\phi(t_{1},\dots,t_{p+1})$ can be computed via construction of a $\Delta$-$\sigma$-Gr\"obner basis of the module $N$ mentioned in the last theorem. Such a basis can be obtained with the use of the generalized Buchberger algorithm described in \cite{Levin3} (where $\leq_{\sigma}$  is to play the role of $\leq_{1}$ in the theory developed in \cite[Section 3]{Levin3}). Another possible way to construct a $\Delta$-$\sigma$-Gr\"obner basis is to use the relative Gr\"obner basis technique developed in \cite{ZW} and \cite{DW}.

The following example deals with a differential field extension ($\sigma = \emptyset$). It illustrates the method of computation of the multivariate dimension polynomial with the use of the last part of Theorem 6.7 and also shows that multivariate dimension polynomials can carry more invariants of a finitely generated differential ($\Delta$-) field extension than the univariate dimension polynomial considered by E. Kolchin.

\begin{example}  {\em Let $K$ be a differential field with a basic set $\Delta =
\{\delta_{1}, \delta_{2}, \delta_{3}\}$ and let
$L = K\langle \eta \rangle$ be a $\Delta$-field extension with the
defining equation $$\delta_{1}^{a+c}\delta_{2}^{b}\eta +
\delta_{2}^{a+b}\delta_{3}^{c}\eta +
\delta_{1}^{a}\delta_{3}^{b+c}\eta = 0$$ where  $a$, $b$, and $c$
are positive integers. (It means that the defining $\Delta$-ideal
of $\eta$ over $K$ is the $\Delta$-ideal of $K\{y\}$ generated
by the $\Delta$-polynomial
$f = \delta_{1}^{a+c}\delta_{2}^{b}y + \delta_{2}^{a+b}\delta_{3}^{c}y
+ \delta_{1}^{a}\delta_{3}^{b+c}y$. Since $f$ is linear, the $\Delta$-ideal $[f]$ is prime.)

As it is shown in \cite[Proposition 6.5.5]{KLMP}, the kernel $N$ of the natural epimorphism of the free $D_{L}$-module $E= D_{L}e$ onto the module of K\"ahler differentials $\Omega_{L|K} = D_{L}d\eta$ is generated by the element
 $g = (\delta_{1}^{a+c}\delta_{2}^{b} + \delta_{2}^{a+b}\delta_{3}^{c} + \delta_{1}^{a}\delta_{3}^{b+c})e$, which, obviously, forms a $\Delta$-Gr\"obner basis of $N$.

Let $\Phi_{\eta}(t_{1}, t_{2}, t_{3})$ be the dimension polynomial associated with the
partition $$\Delta = \{\delta_{1}\}\bigcup\{\delta_{2}\}\bigcup \{\delta_{3}\}.$$
Applying the last part of Theorem 6.7 and the method of evaluation of the sets
$U^{(1)}_{r_{1}, r_{2}, r_{3}}$ and $U^{(2)}_{r_{1}, r_{2}, r_{3}}$ described in the proof of Theorem 3.1, we obtain that

\medskip

\noindent$\Phi_{\eta} =
\left[\D{t_{1}+1\choose 1}{t_{2}+1\choose 1}{t_{3}+1\choose 1} - {t_{1}+1-(a+c)\choose 1}
{t_{2}+1-b\choose 1}{t_{3}+1\choose 1}\right]+ [(t_{1}-(a+c) +1)a(t_{3}+1) + (t_{1}-(a+c) +1)(t_{2}-b +1)(b+c) -
(t_{1}-(a+c) +1)a(b+c)] =$

\medskip

\noindent$(b+c)t_{1}t_{2} + (a+b)t_{1}t_{3} +
(a+c)t_{2}t_{3} +$  terms of total degree at most $1$.

\medskip

Note that the corresponding Kolchin differential dimension polynomial, which describes $tr.deg_{K}K(\Lambda(r)\eta)$,
is as follows:

\medskip

$\phi_{\eta|K}(t) = \D{t+3\choose 3} - \D{t+3- (a+b+c)\choose 3} =
\D\frac{a+b+c}{2}\,t^{2} +$ terms of degree at most $1$.

\medskip
(One can easily obtain this expression either with the use of the classical Gr\"obner basis method or by constructing a free resolution for $\Omega_{L|K}$; the corresponding methods and examples can be found, for example,  in \cite[Chapter 9]{KLMP}.)

\medskip

In this case the Kolchin polynomial $\phi_{\eta|K}(t)$ carries just one differential birational invariant $a+b+c$ while
$\Phi_{\eta}(t_{1}, t_{2}, t_{3})$ determines three such invariants, $a+b, a+c$, and $b+c$, that is, $\Phi_{\eta}$
determines all three parameters of the defining equation while $\phi_{\eta}(t)$ gives just the sum of these parameters.}
\end{example}

\section{Conclusion}
We have proved the existence and presented a method of computation of multivariate dimension polynomials associated with difference-differential field extensions and systems of algebraic difference-differential equations. The proposed method, which generalizes the Ritt-Kolchin characteristic set technique, allows one to evaluate the strength of a system of algebraic difference-differential equations in the sense of A. Einstein. We have also
presented a method of computation of multivariate difference-differential polynomials based on a generalized Gr\"obner basis technique that involves several term orderings. Finally, we have found new invariants of a finitely generated difference-differential field extension that are not carried by univariate Kolchin-type dimension polynomials.

\end{document}